\documentclass[11pt]{amsart}
\usepackage{amscd,amssymb}
\usepackage{amsthm,amsmath,amssymb}
\usepackage[matrix,arrow]{xy}

\newtheorem{theorem}[equation]{Theorem}

\newtheorem{proposition}[equation]{Proposition}
\newtheorem{lemma}[equation]{Lemma}
\newtheorem{corollary}[equation]{Corollary}

\theoremstyle{definition}
\newtheorem{example}[equation]{Example}
\newtheorem{definition}[equation]{Definition}

\theoremstyle{remark}
\newtheorem{remark}[equation]{Remark}

\makeatletter\@addtoreset{equation}{section} \makeatother

\author{Ivan Cheltsov}

\title{Fano varieties with many selfmaps}

\address{\begin{tabbing}
School of Mathematics\\
University of Edinburgh\\
Edinburgh EH9 3JZ, UK\\
\\
\texttt{I.Cheltsov@ed.ac.uk}
\end{tabbing}}

\begin{document}

\begin{abstract}
We study global log canonical thresholds on anticanonically
embedded quasismooth weighted Fano threefold hypersurfaces having
terminal quotient singularities to prove the existence of a
K\"ahler--Einstein metric on most of them, and to produce examples
of Fano varieties with infinite discrete groups of birational
automorphisms.
\end{abstract}

\maketitle

\section{Introduction.}
\label{section:intro}

Let $X$ be a Fano variety\footnote{We assume that all varieties
are projective, normal, and defined over $\mathbb{C}$.} of
dimension $n$ that has at most log terminal singularities.

\begin{definition}
\label{definition:threshold} The global log canonical threshold of
the variety $X$ is the number
$$
\mathrm{lct}\big(X\big)=\mathrm{sup}\left\{\lambda\in\mathbb{Q}\ \left|%
\aligned
&\mathrm{the\ log\ pair}\ \big(X, \lambda H\big)\ \mathrm{has\ log\ canonical\ singularities}\\
&\mathrm{for\ every\ effective}\ \text{$\mathbb{Q}$-divisor}\ H\ \mathrm{such\ that}\ H\equiv -K_{X}\\
\endaligned\right.\right\}\geqslant 0.%
$$
\end{definition}

It follows from \cite{Ti87}, \cite{Na90}, \cite{DeKo01} that the
Fano variety $X$ has an orbifold K\"ahler--Einstein metric in the
case when $X$ has quotient singularities and the inequality
$\mathrm{lct}(X)>n/(n+1)$ holds\footnote{The number
$\mathrm{lct}(X)$ is an algebraic counterpart of the
$\alpha$-invariant introduced in \cite{Ti87}.}.

\begin{example}
\label{example:sextic-double-solid} Let $X$ be a general
hypersurface in $\mathbb{P}(1^{4},3)$ of degree~$6$. Then
$\mathrm{lct}(X)=1$~by~\cite{Pu04d}.
\end{example}

Quasismooth anticanonically embedded weighted Fano threefold
hypersurfaces with terminal singularities  are studied
extensively in \cite{CPR}, \cite{ChPa05}, \cite{Ch06e},
\cite{ChPa05h}. In this paper we prove the following~result.

\begin{theorem}
\label{theorem:main} Let $X$ be a general quasismooth hypersurface
in $\mathbb{P}(1,a_{1},\ldots,a_{4})$ of
degree~$\sum_{i=1}^{4}a_{i}$ having at most terminal singularities
such that $-K_{X}^{3}\leqslant 1$. Then $\mathrm{lct}(X)=1$.
\end{theorem}

The proof of Theorem~\ref{theorem:main} is algebro-geometric, but
Theorem~\ref{theorem:main} implies the following result.

\begin{corollary}
\label{corollary:KE}
With\,the\,assumptions\,of\,Theorem\,\ref{theorem:main},\,the\,variety\,$X$\,has\,a\,K\"ahler--Einstein\,metric.
\end{corollary}

It follows from \cite{CPR}, \cite{Pu04d} that
Theorem~\ref{theorem:main} also implies the following result (see
Theorem~\ref{theorem:Cheltsov}).

\begin{corollary}
\label{corollary:Cheltsov} Let $X_{1},\ldots,X_{r}$ be varieties
that satisfy all hypotheses of Theorem~\ref{theorem:main}.~Then
$$
\mathrm{Bir}\Big(X_{1}\times\cdots\times X_{r}\Big)=\Big<\prod_{i=1}^{r}\mathrm{Bir}\big(X_{i}\big),\ \mathrm{Aut}\Big(X_{1}\times\cdots\times X_{r}\Big)\Big>,%
$$%
the variety $X_{1}\times\cdots\times X_{r}$ is non-rational, and
for any dominant map $\rho\colon X_{1}\times\cdots\times
X_{r}\dasharrow Y$ whose general fiber is rationally connected,
there is a commutative diagram
$$
\xymatrix{
X_{1}\times\cdots\times X_{r}\ar@{->}[d]_{\pi}\ar@{-->}[rr]^{\sigma}&&X_{1}\times\cdots\times X_{r}\ar@{-->}[d]^{\rho}\\
X_{i_{1}}\times\cdots\times X_{i_{k}}\ar@{-->}[rr]_{\xi}&&Y,}%
$$
where $\xi$ and $\sigma$ are
birational maps, and $\pi$ is a  projection for some $\{i_{1},\ldots,i_{k}\}\subsetneq\{1,\ldots,r\}$.%
\end{corollary}

Unlike those of dimension three, no Fano varieties of dimension
four or higher having infinite~groups of birational automorphisms
whose birational automorphisms are well understood have been known
so far. However, we can now easily obtain the following example.

\begin{example}
\label{example:41-41} Let $X$ be a general hypersurface in
$\mathbb{P}(1,1,4,5,10)$ of degree~$20$. Then it immediately
follows from \cite{CPR}, \cite{ChPa05} and
Corollary~\ref{corollary:Cheltsov} that there is an exact sequence
of groups
$$
1\longrightarrow\prod_{i=1}^{m}\Big(\mathbb{Z}_{2}\ast\mathbb{Z}_{2}\Big)\longrightarrow\mathrm{Bir}\Big(\underbrace{X\times\cdots\times X}_{m\ \mathrm{times}}\Big)\longrightarrow\mathrm{S}_{m}\longrightarrow 1,%
$$
where $\mathbb{Z}_{2}\ast\mathbb{Z}_{2}$ is the infinite dihedral
group.
\end{example}

The assertion of Theorem~\ref{theorem:main} may fail without
the~generality~assumption.

\begin{example}
\label{example:34} Let $X$ be a hypersurface in
$\mathbb{P}(1,1,2,6,9)$ of degree $18$ given by the equation
$$
w^{2}=t^{3}+z^{9}+y^{18}+x^{18}\subset\mathbb{P}\big(1,1,2,6,9\big)\cong\mathrm{Proj}\Big(\mathbb{C}[x,y,z,t,w]\Big),
$$
where $\mathrm{wt}(x)=\mathrm{wt}(y)=1$, $\mathrm{wt}(z)=2$,
$\mathrm{wt}(t)=6$, $\mathrm{wt}(w)=9$. The hypersurface $X$ has
terminal quotient singularities, and $-K_{X}^{3}=1/6$.
Arguing~as~in the proof of Theorem~\ref{theorem:main}, we see~that
$$
\mathrm{lct}\big(X\big)=\mathrm{sup}\Big\{\lambda\in\mathbb{Q}\
\Big\vert\ \mathrm{the\ log\ pair}\ \big(X, \lambda D\big)\ \mathrm{is\ log\ canonical\ for\ every\ Weil\ divisor}\ D\in\big|-K_{X}\big|\Big\},%
$$
which easily implies that $\mathrm{lct}(X)=17/18$ by  Lemma~8.12
and Proposition~8.14 in \cite{Ko97}.
\end{example}

Nevertheless, the proof of Theorem~1.3 in \cite{CPR} and the proof
of Theorem~\ref{theorem:main} can also be used to~construct
explicit examples of Fano threefolds to which
Corollaries~\ref{corollary:KE} and \ref{corollary:Cheltsov} can be
applied.

\begin{example}
\label{example:22} Let $X$ be a hypersurface in
$\mathbb{P}(1,2,2,3,7)$ of degree $14$ given by the equation
$$
w^{2}=t^{4}z+y^{7}-z^{7}+x^{14}\subset\mathbb{P}\big(1,2,2,3,7\big)\cong\mathrm{Proj}\Big(\mathbb{C}[x,y,z,t,w]\Big),
$$
where $\mathrm{wt}(x)=1$, $\mathrm{wt}(y)=\mathrm{wt}(z)=2$,
$\mathrm{wt}(t)=3$, $\mathrm{wt}(w)=7$. The hypersurface $X$ has
terminal quotient singularities, and $-K_{X}^{3}=1/6$.
Arguing~as~in the proof of Theorem~\ref{theorem:main}, we see~that
$$
\mathrm{lct}\big(X\big)=\mathrm{sup}\Big\{\lambda\in\mathbb{Q}\
\Big\vert\ \mathrm{the\ log\ pair}\ \big(X, \lambda D\big)\ \mathrm{is\ log\ canonical}\Big\},%
$$
where $D$ is the unique Weil divisor in $|-K_{X}|$. Then
$\mathrm{lct}(X)=1$ by Lemma~8.12 and Proposition~8.14 in
\cite{Ko97}. The threefold $X$~has a~K\"ahler--Einstein metric,
and the group $\mathrm{Bir}(X\times X)$~is~finite.
\end{example}

The proof of Theorem~\ref{theorem:main} is based on the results
obtained in \cite{CPR}, \cite{ChPa05}, \cite{Ch06e},
\cite{ChPa05h}, but it is lengthy, because the hypotheses of
Theorem~\ref{theorem:main} are satisfied for general members of
$90$ out of $95$ familes of quasismooth terminal anticanonically
embedded weighted Fano threefold hypersurfaces (see~\cite{IF00}).

For the convenience of the reader, we organize this paper in the
following way:
\begin{itemize}
\item we prove Theorem~\ref{theorem:main} in Section~\ref{section:the-proof} omitting the proofs of Lemmas~\ref{lemma:smooth-points}, \ref{lemma:singular-points}, \ref{lemma:quadratic-involutions};%
\item we prove auxiliary technical Lemmas~\ref{lemma:smooth-points}, \ref{lemma:singular-points}, \ref{lemma:quadratic-involutions} in Sections~\ref{section:smooth-points}, \ref{section:singular-points}, \ref{section:quadratics-involutions}, respectively;%
\item we consider one important generalization of Corollary~\ref{corollary:Cheltsov} in Section~\ref{section:conic-bundles}.%
\end{itemize}

The author would like to thank J.\,Howie, J.\,Koll\'ar,
L.\,O'Carroll, J.\,Park, A.\,Pukh\-likov,~V.\,Sho\-ku\-rov and the
referees for useful comments. The author is grateful to the IHES
for hospitality.

\section{The proof of main result.}
\label{section:the-proof}

Let $X$ be a general quasismooth hypersurface in
$\mathbb{P}(1,a_{1},a_{2},a_{3},a_{4})$ of degree
$d=\sum_{i=1}^{4}a_{i}$~with terminal singularities, and let
$\gimel\in\{1,\ldots,95\}$ be the ordinal number of
$(a_{1},a_{2},a_{3},a_{4})$ in~the~notation of Table~5 in
\cite{IF00}, where $a_{1}\leqslant a_{2}\leqslant a_{3}\leqslant
a_{4}$. Then $-K_{X}^{3}\leqslant 1\iff\gimel\geqslant 6$.

We suppose that $\gimel\geqslant 6$, but there is $D\in|-nK_{X}|$
such that $(X,\frac{1}{n}D)$ is not log canonical, where $n$ is a
natural number. Then to prove Theorem~\ref{theorem:main} it is
enough to derive~a~contradiction, because the class group of the
hypersurface $X$ is generated by the divisor $-K_{X}$.

\begin{remark}
\label{remark:weighted-sum} Let $V$ be a variety, let $B$ and
$B^{\prime}$ be effective $\mathbb{Q}$-Cartier
$\mathbb{Q}$-divisors on  $V$ such that the~singularities of the
log pairs $(V, B)$ and $(V, B^{\prime})$ are log canonical, and
let $\alpha$ be a rational number such that
$0\leqslant\alpha\leqslant 1$. Then the log pair $(V,\ \alpha
B+(1-\alpha)B^{\prime})$ is log canonical.
\end{remark}

Thus, we may assume that $D$ is an irreducible surface due to
Remark~\ref{remark:weighted-sum}.

\begin{lemma}
\label{lemma:anticanclass} The inequality $n\ne 1$ holds.
\end{lemma}

\begin{proof}
Suppose that $n=1$. Then the log pair $(X, D)$ is log canonical at
every singular point of the threefold $X$ by Lemma~8.12 and
Proposition~8.14 in \cite{Ko97}. Thus, the equality $a_{1}=1$
holds, because the linear system $|-K_{X}|$ consists of a single
surface in the case when $a_{1}\ne 1$.

The equality $a_{1}=1$ holds for $36$ values of
$\gimel\in\{6,7,\ldots,95\}$, but all possible cases are very
similar. So for the sake of simplicity, we assume that
$\gimel=14$. Then there is a natural double cover $\pi\colon
X\to\mathbb{P}(1,1,1,4)$~branched~over a general hypersurface
$F\subset \mathbb{P}(1,1,1,4)$ of degree $12$.

Suppose that the singularities of the log pair $(X, D)$ are not
log canonical at~some smooth point $P$ of the threefold $X$. Let
us show that this assumption leads to a contradiction.

Put $\bar{D}=\pi(D)$ and $\bar{P}=\pi(P)$. Counting parameters, we
see that $\mathrm{mult}_{\bar{P}}(F\vert_{\bar{D}})\leqslant 2$,
which is a contradiction, because
$(\bar{D},\frac{1}{2}F\vert_{\bar{D}})$ is not log~ca\-no\-ni\-cal
at $\bar{P}$ by Lemma~8.12~in~\cite{Ko97}.
\end{proof}

\begin{lemma}
\label{lemma:smooth-points} The log pair $(X,\frac{1}{n}D)$ is log
canonical at smooth points of the threefold $X$.
\end{lemma}

\begin{proof}
See Section~\ref{section:smooth-points}.
\end{proof}

Therefore, there is a singular point $O$ of the threefold $X$ such
that $(X,\frac{1}{n}D)$ is not canonical at the point $O$. It
follows from \cite{IF00} that $O$ is a singular point of type
$\frac{1}{r}(1,a,r-a)$, where $a$ and~$r$ are coprime natural
numbers such that $r>2a$ (see Table~5~in~\cite{IF00} for the
values of $a$ and $r$).

Let $\alpha\colon U\to X$ be a blow up of $O$ with weights
$(1,a,r-a)$.~Then
\begin{equation}
\label{equation:degree-of-blow-up}
-K_{U}^{3}=-K_{X}^{3}-\frac{1}{r^{3}}E^{3}=-K_{X}^{3}-\frac{1}{ra(r-a)}=\frac{\sum_{i=1}^{4}a_{i}}{a_{1}a_{2}a_{3}a_{4}}-\frac{1}{ra(r-a)},
\end{equation}
where $E$ is the exceptional divisor of $\alpha$. There is a
rational number $\mu$ such that
$$
\bar{D}\equiv  \alpha^{*}\big(D\big)-\mu E\equiv -nK_{U}+\big(n/r-\mu\big)E
$$
where $\bar{D}$ is the proper transform of $D$ on $U$. Then it
follows from \cite{Ka96} that $\mu>n/r$.

\begin{lemma}
\label{lemma:negative-K-cube} The inequality $-K_{U}^{3}\geqslant 0$ holds.%
\end{lemma}

\begin{proof}
Suppose that $-K_{U}^{3}<0$. Let $C$ be a curve in $E$. Then the
curve $C$ generates an extremal ray of the cone $\mathbb{NE}(U)$.
Moreover, it follows from Corollary~5.4.6 in \cite{CPR} that there
is an~irreducible curve $\Gamma\subset U$ such that $\Gamma$
generates the extremal ray of $\mathbb{NE}(U)$ that is different
from $\mathbb{R}_{\geqslant 0}C$, and
$$
\Gamma\equiv-K_{U}\cdot \Big(-bK_{U}+cE\Big),
$$
where $b>0$ and $c\geqslant 0$ are integers (see Remark~5.4.7 in
\cite{CPR}).

Let $T$ be a divisor in $|-K_{U}|$. Then $\bar{D}\cdot T$ is
effective, because $\bar{D}\ne T$. However we have
$$
\bar{D}\cdot T\equiv -K_{U}\cdot\Big(-nK_{U}+\big(n/r-\mu\big)E\Big)\not\in\mathbb{NE}(U),%
$$
because $\mu>n/r$, $b>0$, and $c\geqslant 0$. So we have a
contradiction.
\end{proof}

Taking into account the possible values of
$(a_{1},a_{2},a_{3},a_{4})$, we see that
$\gimel\not\in\{75,84,87,93\}$.

\begin{lemma}
\label{lemma:zero-K-cube} The inequality $-K_{U}^{3}\ne 0$ holds.
\end{lemma}

\begin{proof}
Firstly, suppose that $-K_{U}^{3}=0$ and $\gimel\ne 82$. Then the
linear system $|-rK_{U}|$ does not have base points for $r\gg 0$
and induces a morphism $\eta\colon U\to\mathbb{P}(1,a_{1},a_{2})$
such that the diagram
$$
\xymatrix{
&&&U\ar@{->}[lld]_{\alpha}\ar@{->}[rrd]^{\eta}&&&\\%
&X\ar@{-->}[rrrr]_{\psi}&&&&\mathbb{P}(1,a_{1},a_{2})&}
$$
is commutative, where $\psi$ is a natural projection. The morphism
$\eta$ is an elliptic fibration. Thus
$$
\bar{D}\cdot C=-nK_{U}\cdot C+\big(n/r-\mu\big)E\cdot C=\big(n/r-\mu\big)E\cdot C<0,%
$$
where $C$ is a general fiber of $\eta$, which is a contradiction.

Suppose that  $-K_{U}^{3}=0$ and $\gimel=82$. Then $X$ is a
hypersurface in $\mathbb{P}(1,1,5,12,18)$ of degree $36$, whose
singularities consist of two points $P$ and $Q$ of types
$\frac{1}{5}(1,2,3)$ and $\frac{1}{6}(1,1,5)$, respectively.

We see that either $P=O$, or $Q=O$. The hypersurface $X$ can be
given by the equation
$$
z^{7}y+\sum_{i=0}^{6}z^{i}f_{36-5i}\big(x,y,z,t\big)=0\subset\mathbb{P}\big(1,1,5,12,18\big)\cong\mathrm{Proj}\Big(\mathbb{C}[x,y,z,t,w]\Big),
$$
where $\mathrm{wt}(x)=\mathrm{wt}(y)=1$, $\mathrm{wt}(z)=5$,
$\mathrm{wt}(t)=12$, $\mathrm{wt}(w)=18$, and $f_{i}$ is a
quasi\-ho\-mo\-ge\-ne\-ous polynomial of degree $i$. Then $P$ is
is given by the equations $x=z=t=w=0$.

Suppose that $Q=O$. Then  the linear system $|-rK_{U}|$ has no
base points for $r\gg 0$, which leads to a contradiction as in the
case when $\gimel\ne 82$. So we see that $P=O$.

Let $\bar{S}$ be the proper transform on $U$ of the surface that
is cut out on $X$ by $y=0$. Then
$$
\bar{S}\equiv  \alpha^{*}\big(-K_{X}\big)-\frac{6}{5}E,
$$
and the base locus of the pencil $|-K_{U}|$ consists of two
irreducible curves $L$ and $C$ such that the~curve $L$ is
contained in the $\alpha$-exceptional surface $E$, and the curve
$\pi(C)$ is the unique base curve of the~pencil $|-K_{X}|$. Then
$-K_{U}\cdot C=-1/6$ and $-K_{U}\cdot L>0$. We have $\mu\leqslant
n/5$ due to
$$
n/5-\mu=\Big(-K_{U}+\alpha^{*}\big(-5K_{X}\big)\Big)\cdot\bar{S}\cdot\bar{D}\geqslant 0,%
$$
because it follows from Lemma~8.12 and Proposition~8.14 in
\cite{Ko97} that $\bar{D}\ne\bar{S}$. However we  know that the
inequality $\mu>n/5$ holds by \cite{Ka96}. So again we have a
contradiction.
\end{proof}

Thus, taking into account the
equality~\ref{equation:degree-of-blow-up} and possible values of
$(a_{1},a_{2},a_{3},a_{4})$, we see that
$$
\gimel\not\in\Big\{11, 14, 19, 22, 28, 34, 37, 39, 49, 52, 53, 57, 59, 64, 66, 70, 72, 73, 78, 80, 81, 86, 88, 89, 90, 92, 94, 95\Big\}%
$$
by Lemma~\ref{lemma:zero-K-cube}. So the assertion of
Theorem~\ref{theorem:main} is proved for $32$ values
of~$(a_{1},a_{2},a_{3},a_{4})$.

\begin{lemma}
\label{lemma:superrigid} The groups $\mathrm{Bir}(X)$ and
$\mathrm{Aut}(X)$ do not coincide.
\end{lemma}

\begin{proof}
Suppose that $\mathrm{Bir}(X)=\mathrm{Aut}(X)$. Let $\bar{S}$ be a
general surface in  $|-K_{U}|$. Then it follows from Lemma~5.4.5
in \cite{CPR} that there is an irreducible surface $\bar{T}\subset
U$ such that
\begin{itemize}
\item the equivalence $\bar{T}\sim c\bar{S}-bE$ holds, where $c\geqslant 1$ and $b\geqslant 1$ are natural numbers,%
\item the scheme-theoretic intersection $\bar{T}\cdot\bar{S}$ is an irreducible and reduced curve $\Gamma$,%
\item the curve $\Gamma$ generates an extremal ray of the cone $\mathbb{NE}(U)$.%
\end{itemize}

The surface $\bar{T}$ is easy to construct explicitly (see
\cite{CPR}), and the possible values for the natural numbers $c$
and $b$ can be found in \cite{CPR}. The surface $\bar{T}$ is
determined uniquely by the point $O$.

Put $T=\alpha(\bar{T})$. Then it follows from Lemma~8.12 and
Proposition~8.14 in \cite{Ko97} that~the~singularities of the log
pair $(X,\frac{1}{c}T)$ are log canonical. Therefore, we have
$D\ne T$.

Let $\mathcal{P}$ be the pencil~gene\-rated by the effective
divisors $nT$ and $cD$. Then the singularities of the log pair
$(X, \frac{1}{cn}\mathcal{P})$~are~not~canonical, which is
impossible due to \cite{CPR}.
\end{proof}

It follows from \cite{CPR} and Lemma~\ref{lemma:superrigid} that
$\gimel\not\in\{11, 21, 29, 35, 50, 51, 55, 62, 63, 67, 71, 77,
82, 83, 85, 91\}$.

\begin{lemma}
\label{lemma:nef-and-big} The divisor $-K_{U}$ is nef.
\end{lemma}

\begin{proof}
Suppose that $-K_{U}$ is not nef. Then it follows from \cite{CPR}
that $\gimel=47$ and $O$ is a singular point of type
$\frac{1}{5}(1,2,3)$. The hypersurface $X$ can be given by the
equation
$$
z^{4}y+\sum_{i=0}^{3}z^{i}f_{21-5i}\big(x,y,z,t\big)=0\subset\mathbb{P}\big(1,1,5,7,8\big)\cong\mathrm{Proj}\Big(\mathbb{C}[x,y,z,t,w]\Big),
$$
where $\mathrm{wt}(x)=\mathrm{wt}(y)=1$, $\mathrm{wt}(z)=5$,
$\mathrm{wt}(t)=7$, $\mathrm{wt}(w)=8$, and $f_{i}$ is a general
quasi\-ho\-mo\-ge\-ne\-ous polynomial of degree $i$. Let $S$ be
the surface on $X$ that is cut out by the equation $y=0$, and
$\bar{S}$ be the proper transform of the surface $S$ on the
threefold $U$. Then
$$
\bar{S}\equiv  \alpha^{*}\big(-K_{X}\big)-\frac{6}{5}E,
$$
but the divisor $-3K_{U}+\alpha^{*}(-5K_{X})$ is nef (see
\cite{Ch06e}). Thus, the inequality $\mu\leqslant n/5$ holds due
to
$$
n/5-\mu=\frac{1}{3}\Big(-3K_{U}+\alpha^{*}\big(-5K_{X}\big)\Big)\cdot\bar{S}\cdot\bar{D}\geqslant 0,%
$$
because $D\ne S$. However we know that $\mu>n/5$. So we have a
contradiction.
\end{proof}

Thus, the divisor $-K_{U}$ is nef and big, because $-K_{U}^{3}>0$
by Lemmas~\ref{lemma:negative-K-cube} and \ref{lemma:zero-K-cube}.

\begin{lemma}
\label{lemma:discrepancy} The inequality $\mu/n-1/r<1$ holds.
\end{lemma}

\begin{proof}
We only consider the case when $\gimel=58$ and $O$ is a
sin\-gu\-lar point of type $\frac{1}{10}(1,3,7)$, because the
proof is similar in all other cases (cf.
Lemma~\ref{lemma:inequality-for-mu}). Then $X$ can be given~by
$$
w^{2}z+wf_{14}\big(x,y,z,t\big)+f_{24}\big(x,y,z,t\big)=0\subset\mathbb{P}\big(1,3,4,7,10\big)\cong\mathrm{Proj}\Big(\mathbb{C}[x,y,z,t,w]\Big),
$$
where $\mathrm{wt}(x)=1$, $\mathrm{wt}(y)=3$, $\mathrm{wt}(z)=4$,
$\mathrm{wt}(t)=7$, $\mathrm{wt}(w)=10$, and $f_{i}$ is a
quasi\-ho\-mo\-ge\-ne\-ous polynomial of degree $i$. Let $R$ be
the surface on $X$ that is cut out by  $t=0$, and $\bar{R}$ be the
proper transform of the surface $R$ on the threefold $U$. Then
$$
\bar{R}\equiv  \alpha^{*}\big(-4K_{X}\big)-\frac{7}{5}E,
$$
and $(X, \frac{1}{4}R)$ is log canonical at $O$ by Lemma~8.12 and
Proposition~8.14 in \cite{Ko97}. Then $R\ne D$~and
$$
0\leqslant -K_{U}\cdot\bar{R}\cdot\bar{D}=4n/35-2\mu/3,%
$$
because $-K_{U}$ is nef. Thus, we have $\mu\leqslant 6n/35$, which
implies that $\mu/n-1/10<1$.
\end{proof}

So the log pair $(U, \frac{1}{n}\bar{D}+(\mu/n-1/r)E)$ is not log
canonical at some point $P\in E$,~because
$$
K_{U}+\frac{1}{n}\bar{D}\equiv \alpha^{*}\Big(K_{X}+\frac{1}{n}D\Big)+\big(1/r-\mu/n\big)E.%
$$

\begin{lemma}
\label{lemma:singular-points} The threefold $U$ is smooth at the
point $P$.
\end{lemma}

\begin{proof}
See Section~\ref{section:singular-points}.
\end{proof}

Thus, the inequality $\mathrm{mult}_{P}(\bar{D})>n+n/r-\mu$ holds.
But it follows from \cite{CPR} that
\begin{itemize}
\item either $d=2r+a_{j}$ for some $j$, and there is a quadratic involution $\tau\in\mathrm{Bir}(X)$ induced~by~$O$,%
\item or $d=3r+a_{j}$ for some $j$, and there is an elliptic involution $\tau\in\mathrm{Bir}(X)$ induced by $O$,%
\end{itemize}
where $d=\sum_{i=1}^{4}a_{i}$.

\begin{lemma}
\label{lemma:quadratic-involutions} The inequality $d\ne 2r+a_{j}$
holds for every $j\in\{1,2,3,4\}$.
\end{lemma}

\begin{proof}
See Section~\ref{section:quadratics-involutions}.
\end{proof}

Thus, it follows from \cite{CPR} that there is $j\in\{1,2,3,4\}$
such that $d=3r+a_{j}$.

\begin{remark}
\label{remark:curves} Let $V$ be a threefold with isolated
singularities, let $B\ne T$ be effective irreducible divisors on
the threefold $V$, and let $H$~be a nef divisor on the threefold
$V$. Put
$$
B\cdot T=\sum_{i=1}^{r}\epsilon_{i} L_{i}+\Delta,
$$
where $L_{i}$ is an irreducible curve, $\epsilon_{i}$ is a
non-negative integer, and $\Delta$ is an effective one-cycle whose
support does not contain the curves $L_{1},\ldots,L_{r}$. Then
$\sum_{i=1}^{r}\epsilon_{i} H\cdot L_{i}\leqslant B\cdot T\cdot
H$.
\end{remark}

It follows from Lemma~\ref{lemma:quadratic-involutions} that
$\gimel\in\{7, 20, 23, 36, 40, 44, 61, 76\}$ (see \cite{CPR}).

\begin{lemma}
\label{lemma:7-20-36-elliptic} The case $\gimel\in\{7,20,36\}$ is
impossible.
\end{lemma}

\begin{proof}
Suppose that $\gimel\in\{7,20,36\}$. Then $a_{1}=1$, and it
follows from Lemma~\ref{lemma:quadratic-involutions} that
$O$~is~a~singular point of type $\frac{1}{a_{2}}(1,1,a_{2}-1)$.
Then $|-rK_{U}|$ induces a birational morphism $\sigma\colon U\to
V$ such that $\sigma$ contracts smooth rational curves
$C_{1},\ldots, C_{l}$, and $V$ is a hypersurface in
$\mathbb{P}(1,1,a_{3},2a_{4},3a_{4})$ of degree $6a_{4}$, where
$l=d(d-a_{4})/a_{3}$. Let $T$ be the surface in $|-K_{U}|$ that
contains $P$.

Suppose that $P\not\in\cup_{i=1}^{l}C_{i}$. Then it follows from
the proof of Theorem~5.6.2 in \cite{CPR} that there are natural
number $s>0$ and a surface $H\in|-s2a_{4}K_{U}|$ such that
$$
s2a_{4}\big(-nK_{X}^{3}-\mu/(a_{2}-1)\big)=\bar{D}\cdot T\cdot H\geqslant \mathrm{mult}_{P}\big(\bar{D}\big)s>s\big(n+n/a_{2}-\mu\big) s,%
$$
which is impossible, because $\mu>n/a_{2}$. So we may assume that
$P\in C_{1}$.

Put $\bar{D}\cdot T=mC_{1}+\Delta$, where $m$ is a non-negative
integer, and $\Delta$ is an effective cycle such that the support
of $\Delta$ does not~contain the~curve $C_{1}$. The curve $C_{1}$
is a smooth rational curve such that $\alpha^{*}(-K_{X})\cdot
C_{1}=2/a_{2}$ and $E\cdot C_{1}=2$.

It follows from \cite{CPR} that there is a surface
$R\in|-a_{3}K_{U}|$ such that~$R$~contains $C_{1}$, but $R$~does
not contain components of the cycle $\Delta$ passing through the
point $P$. Then
$$
a_{3}\big(-nK_{X}^{3}-\mu/(a_{2}-1)\big)=R\cdot\Delta\geqslant \mathrm{mult}_{P}\big(\Delta\big)>n+n/a_{2}-\mu-m,%
$$
which implies that $m>a_{3}n/a_{4}$, because $\mu>n/a_{2}$.
Therefore, we have
$$
a_{3}n/a_{4}<m\leqslant\frac{-dnK_{X}\cdot\alpha(C_{1})}{a_{1}a_{2}a_{3}a_{4}}=\frac{dn}{2a_{1}a_{3}a_{4}}%
$$
by Remark~\ref{remark:curves}, because
$-K_{X}\cdot\alpha(C_{1})=2/a_{2}$. The inequalities just obtained
imply that $\gimel=7$.

Let $\psi\colon X\dasharrow\mathbb{P}(1,a_{1},a_{2})$ be a natural
projection. The fiber of $\psi$ over  $\psi(P)$ consists of two
irreducible components, and one of them is $C_{1}$. Let $Z$ be the
other component of this fiber.~Then
$$
C_{1}^{2}=-2,\ C_{1}\cdot Z=2,\ Z^{2}=-4/3
$$
on the surface $T$. Put $\Delta=\bar{m}Z+\Omega$, where $\bar{m}$
is a non-negative integer and $\Omega$~is~an~effective one-cycle
whose support does not~contain the~curve $Z$. Then
$$
4n/3-2\mu-5\bar{m}/3=\big(Z+C_{1}\big)\cdot\Omega>3n/2-\mu-m,
$$
and $4\bar{m}/3\geqslant 2m-5n/6$, because $\Omega\cdot Z\geqslant
0$. The inequalities just obtained immediately  imply that the
inequality $\mu\leqslant n/2$ holds. So we have a contradiction,
because $\mu>n/2$.
\end{proof}

Hence, it follows from Lemmas~\ref{lemma:quadratic-involutions}
and \ref{lemma:7-20-36-elliptic}~that
$\gimel\in\{23,40,44,61,76\}$ and~$d=3r+a_{j}$, where $r=a_{3}>2a$
and $1\leqslant j\leqslant 2$. Then $X$ has~a~singular point $Q$
of type $\frac{1}{\bar{r}}(1,\bar{a},\bar{r}-\bar{a})$~such~that
$$
-K_{X}^{3}=\frac{1}{ra(r-a)}+\frac{1}{\bar{r}\bar{a}(\bar{r}-\bar{a})},%
$$
where $\bar{r}=a_{4}>2\bar{a}$ and $\bar{a}\in\mathbb{N}$. It
follows from \cite{CPR} that there is a commutative diagram
$$
\xymatrix{
&U\ar@{->}[ddl]_{\sigma}\ar@{->}[rd]^{\alpha}&&&W\ar@{->}[lll]_{\gamma}\ar@{->}[rrdd]^{\eta}&\\
&&X\ar@{-->}[drrrr]^{\psi}&&\\
V\ar@{^{(}->}[rr]&&\mathbb{P}\big(1,a_{1},a_{2},2a_{4},3a_{4}\big)\ar@{-->}[rr]_{\chi}&&\mathbb{P}\big(1,a_{1},a_{2},2a_{4}\big)\ar@{-->}[rr]_{\xi}&&\mathbb{P}\big(1,a_{1},a_{2}\big),}
$$
where $\xi$, $\chi$, $\psi$ are projections, $\eta$ is an elliptic
fibration, $\gamma$ is a weighted blow up of a point that
dominates the point $Q$ with weights $(1,\bar{a},
\bar{r}-\bar{a})$, and $\sigma$ is a birational morphism that
contracts smooth curves $C_{1},\ldots, C_{l}$ such that $V$ is a
hypersurface in $\mathbb{P}(1,a_{1},a_{2},2a_{4},3a_{4})$ of
degree $6a_{4}$, where $l=d(d-a_{4})/(a_{1}a_{2})$. Let $L$ be a
curve in $|\mathcal{O}_{\mathbb{P}(1,\,a,\,r-a)}(1)|$, where
$E\cong \mathbb{P}(1,a,r-a)$.

\begin{lemma}
\label{lemma:lower-bound-for-mu} Suppose that $P\not\in L$. Then
$\mu>na(r+1)/(r^2+ar)$.
\end{lemma}

\begin{proof}
There is a~unique curve
$C\in|\mathcal{O}_{\mathbb{P}(1,\,a,\,r-a)}(a)|$ such that $P\in
C$. Put $\bar{D}\big\vert_{E}=\delta C+\Upsilon\equiv r\mu L$,
where $\delta$ is a non-negative integer and $\Upsilon$ is an
effective cycle such that $C\not\subset\mathrm{Supp}(\Upsilon)$.
Then
$$
(r\mu-a\delta)/(r-a)=\big(r\mu-a\delta\big)L\cdot C=C\cdot\Upsilon\geqslant\mathrm{mult}_{P}\big(\Upsilon\big)>n+n/r-\mu-\delta,%
$$
which implies that $\mu>na(r+1)/(r^2+ar)$, because
$\delta\leqslant r\mu/a$.
\end{proof}

Let $T$ be a surface in $|-K_{U}|$. Then $-K_{U}\cdot T
\cdot\bar{D}\geqslant 0$, which implies that $\mu\leqslant
-na(r-a)K_{X}^{3}$.

\begin{lemma}
\label{lemma:elliptic-involution-T} The point $P$ is not contained
in the surface $T$.
\end{lemma}

\begin{proof}
Suppose that $P$ is contained in the surface $T$. Then $P$ is not
contained in the base locus of the pencil $|-a_{1}K_{U}|$, because
the base locus of the pencil $|-a_{1}K_{U}|$ does not contain
smooth points of the surface $E$. The point $P$ is not contained
in the union $\cup_{i=1}^{l}C_{i}$, because $P\in T$.

The proof of Theorem~5.6.2 in \cite{CPR} implies the existence of
a~surface $H\in|-s2a_{1}a_{4}K_{U}|$~such~that
$$
s2a_{1}a_{4}\big(-nK_{X}^{3}-\mu/a_{2}\big)=\bar{D}\cdot H\cdot T\geqslant \mathrm{mult}_{P}\big(\bar{D}\big)s>s\big(n+n/r-\mu\big) s,%
$$
where $s$ is a natural number, which is impossible, because
$\mu>n/r$.
\end{proof}

We have $T\vert_{E}\sim \mathcal{O}_{\mathbb{P}(1,\,a,\,r-a)}(1)$.
Thus, taking into account that $\gimel\in\{7, 20, 23, 36, 40, 44,
61, 76\}$, we see that $\gimel\in\{23, 44\}$ by
Lemmas~\ref{lemma:lower-bound-for-mu} and
\ref{lemma:elliptic-involution-T}, because $\mu\leqslant
-na(r-a)K_{X}^{3}$.

Let $S$ be a surface in $|-a_{1}K_{U}|$ that contains $P$. Then
$\bar{D}\ne S$, because $\mu>n/r$.

\begin{lemma}
\label{lemma:elliptic-involution-exceptional-curves} The point $P$
is contained in $\cup_{i=1}^{l}C_{i}$.
\end{lemma}

\begin{proof}
Suppose that $P\not\in\cup_{i=1}^{l}C_{i}$. Then the proof of
Theorem~5.6.2 in \cite{CPR} implies that
$$
s2a_{1}a_{4}\big(-nK_{X}^{3}-\mu/a_{2}\big)=\bar{D}\cdot H\cdot S\geqslant \mathrm{mult}_{P}\big(\bar{D}\big)s>s\big(n+n/r-\mu\big) s%
$$
for some $s\in\mathbb{N}$ and a surface $H\in|-s2a_{4}K_{U}|$,
which is impossible, because $\mu>n/r$.
\end{proof}

We may assume that $P\in C_{1}$. Put $\bar{D}\cdot
S=mC_{1}+\Delta$, where $m$ is a non-negative integer,
and~$\Delta$~is~an effective cycle whose support does not~contain
$C_{1}$. Then it follows from Remark~\ref{remark:curves} that the
inequality $m\leqslant nd/(a_{2}d-a_{2}a_{3})$ holds, because
$-K_{X}\cdot\alpha(C_{1})=(d-a_{3})/(a_{3}a_{4})$.

It follows from the proof of Theorem~5.6.2 in \cite{CPR} that
there is $R\in|-s2a_{4}K_{U}|$~such~that
$$
s2a_{1}a_{4}\big(-nK_{X}^{3}-\mu/a_{2}\big)=R\cdot\Delta\geqslant\mathrm{mult}_{P}\big(\Delta\big)s>s\big(n+n/r-\mu-m\big),%
$$
where $s\in\mathbb{N}$. However we have $m\leqslant
nd/(a_{2}d-a_{2}a_{3})$, which implies that $\gimel=23$.

Therefore, we proved that $X$ is a hypersurface
$\mathbb{P}(1,2,3,4,5)$ of degree $14$ and $O$~is~a~singular point
of type $\frac{1}{4}(1,1,3)$. Let $M$ be a general surface in the
linear system $|-3K_{X}|$ that passes through the point $P$. Then
$S\cdot M=C_{1}+Z_{1}$, where $Z_{1}$ is a curve such that
$-K_{U}\cdot Z_{1}=1/5$.~Put
$$
\bar{D}\cdot S=mC_{1}+\bar{m}Z_{1}+\Upsilon,
$$
where $\bar{m}$ is a non-negative integer, and $\Upsilon$ is an
effective cycle, whose support does not contain the curves $C_{1}$
and $Z_{1}$. Then $m<7n/15$ by Remark~\ref{remark:curves}. But
$\mu>n/4$~and
$$
7n/10-6\mu/3-3\bar{m}/5=M\cdot\Upsilon\geqslant\mathrm{mult}_{P}\big(\Upsilon\big)>5n/4-\mu-m,%
$$
because $P\not\in Z_{1}$. The inequality obtained implies a
contradiction.

Therefore, the assertion of Theorem~\ref{theorem:main} is
completely proved.

\section{Non-singular points.}
\label{section:smooth-points}

In this section we prove the assertion of
Lemma~\ref{lemma:smooth-points}. Let us use the assumptions and
notation of Lemma~\ref{lemma:smooth-points}. Take an arbitrary
smooth point $P$ of  the threefold $X$.

\begin{lemma}
\label{lemma:smooth-points-small-K-cube} Suppose that $a_{4}$
divides $d$, $a_{1}\ne a_{2}$, and $-a_{2}a_{3}K_{X}^{3}\leqslant
1$. Then $\mathrm{mult}_{P}(D)\leqslant n$.
\end{lemma}

\begin{proof}
Suppose that $\mathrm{mult}_{P}(D)>n$. Let $L$ be  the base curve
of $|-a_{1}K_{X}|$, and $T$ be~a~surface in the linear system
$|-K_{X}|$. Then $D\cdot T$ is an effective one-cycle, and
$\mathrm{mult}_{P}(L)=1$.

Suppose that $P\in L$. Let $R$ be a general surface in
$|-a_{1}K_{X}|$. Put $D\cdot T=mL+\Delta$, where $m$ is
non-negative integer, and $\Delta$ is an effective cycle whose
support does not contain $L$. Then
$$
-a_{1}\big(n-a_{1}m\big)K_{X}^{3}=D\cdot T\cdot R-mR\cdot L=R\cdot\Delta\geqslant\mathrm{mult}_{P}\big(\Delta\big)>n-m,%
$$
which is impossible, because $-a_{1}K_{X}^{3}\leqslant 1$. Thus,
we see that $P\not\in L$.

Suppose that $P\in T$. It follows from Theorem~5.6.2 in \cite{CPR}
that
$$
ns\geqslant -sa_{1}a_{3}nK_{X}^{3}=D\cdot S\cdot T\geqslant\mathrm{mult}_{P}\big(D\big)s>ns%
$$
for some $s\in\mathbb{N}$ and some surface $S\in |-sa_1a_3K_{X}|$.
Hence, we see that $P\not\in T$.

Let $G$ be a general surface in $|-a_{2}K_{X}|$ that contains $P$.
Then $G\cdot D$ is an effective~cycle, but it follows from
Theorem~5.6.2 in \cite{CPR} that there are $s\in\mathbb{N}$ and
$H\in|-sa_3K_{X}|$~such that
$$
ns\geqslant-sa_{2}a_{3}nK_{X}^{3}=D\cdot H\cdot G\geqslant \mathrm{mult}_{P}\big(D\big)s>ns,%
$$
because $-a_{2}a_{3}K_{X}^{3}\leqslant 1$. So we have a
contradiction.
\end{proof}

\begin{lemma}
\label{lemma:smooth-points-11-34-53-88} Suppose that $a_{4}$
divides $d$, $1=a_{1}\ne a_{2}$, and $-a_{3}K_{X}^{3}\leqslant 1$.
Then $\mathrm{mult}_{P}(D)\leqslant n$.
\end{lemma}

\begin{proof}
Suppose that $\mathrm{mult}_{P}(D)>n$. Then arguing as in the
proof of Lemma~\ref{lemma:smooth-points-small-K-cube}, we see that
the point $P$ is not contained in the base curve of the pencil
$|-K_{X}|$.

Let $T$ be a general surface in $|-K_{X}|$ that contains $P$. Then
Theorem~5.6.2 in \cite{CPR} implies~that there are
$s\in\mathbb{N}$ and $S\in|-sa_3K_{X}|$ such that
$ns\geqslant-sa_{3}nK_{X}^{3}=D\cdot S\cdot T\geqslant
\mathrm{mult}_{P}\big(D\big)s>ns$.
\end{proof}

\begin{lemma}
\label{lemma:smooth-points-a1-a4} Suppose that $a_{1}\ne a_{2}$
and $-a_{1}a_{4}K_{X}^{3}\leqslant 1$. Then
$\mathrm{mult}_{P}(D)\leqslant n$.
\end{lemma}

\begin{proof}
Suppose that $\mathrm{mult}_{P}(D)>n$.
The~proof~of~Lemma~\ref{lemma:smooth-points-11-34-53-88} implies
that $a_{1}\ne 1$. Arguing as~in the proof of
Lemma~\ref{lemma:smooth-points-small-K-cube}, we see that $P$ is
not contained in the~unique surface of~$|-K_{X}|$.

Let $S$ be a surface in $|-a_{1}K_{X}|$ that contains $P$. Then we
may~assume~that $\mathrm{mult}_{P}\big(S\big)\leqslant a_{1}$,
because $P\not\in T$ and $X$ is sufficiently general. Thus, we
have $S\ne D$.

It follows from Theorem~5.6.2 in \cite{CPR} that there are
$s\in\mathbb{N}$ and $H\in|-sa_4K_{X}|$~such~that~$H$~has
multiplicity at least $s>0$ at $P$ and contains no components of
$D\cdot S$ passing through $P$.~Then
$$
ns\geqslant-sa_{1}a_{4}nK_{X}^{3}=D\cdot S\cdot H\geqslant \mathrm{mult}_{P}\big(D\big)s>ns,%
$$
because $-a_{1}a_{4}K_{X}^{3}\leqslant 1$. So we have a
contradiction.
\end{proof}

Taking into account the possible values of
$(a_{1},a_{2},a_{3},a_{4})$, we see~that
$\mathrm{mult}_{P}(D)\leqslant n$ whenever
$$
\gimel\not\in\Big\{6,7,8,9,10,12,13,14,16,18,19,20,22,23,24,25,32,33,38\Big\}
$$
by Lemmas~\ref{lemma:smooth-points-small-K-cube},
\ref{lemma:smooth-points-11-34-53-88},
\ref{lemma:smooth-points-a1-a4}. The log pair $(X, \frac{1}{n}D)$
is log canonical at $P$ if $\mathrm{mult}_{P}(D)\leqslant n$ (see
\cite{Ko97}).

\begin{lemma}
\label{lemma:smooth-points-18} Suppose that $\gimel=18$. Then $(X,
\frac{1}{n}D)$ is log canonical at $P$.
\end{lemma}

\begin{proof}
Suppose that the log pair $(X, \frac{1}{n}D)$ is log canonical at
the point $P$. Let us show that this assumption leads to a
contradiction. Note that the inequality $\mathrm{mult}_{P}(D)>n$
holds.

The threefold $X$ is a hypersurface in $\mathbb{P}(1,2,2,3,5)$ of
degree $12$, whose~singularities consist of six points of type
$\frac{1}{2}(1,1,1)$, and a point $O$ of type
$\frac{1}{5}(1,2,3)$. It follows from \cite{Ch06e} that the
diagram
$$
\xymatrix{
&U\ar@{->}[d]_{\alpha}&&W\ar@{->}[ll]_{\beta}\ar@{->}[d]^{\eta}&\\%
&X\ar@{-->}[rr]_{\psi}&&\mathbb{P}(1,2,2)&}
$$
commutes, where $\alpha$ is a weighted blow up of the point $O$
with weights $(1,2,3)$, $\beta$ is a weighted blow up with weights
$(1,1,3)$ of a singular point of type $\frac{1}{3}(1,1,2)$, and
$\eta$ is an elliptic fibration.

Let $C$ be a fiber of the projection $\psi$ that passes through
the point $P$, and $L$ be its irreducible reduced component. We
have $-K_{X}\cdot C=4/5$. But the number $-5K_{X}\cdot L$ is
natural if $L$ contains no points of type $\frac{1}{2}(1,1,1)$.
Then $C=2L$ whenever $C$ contains a point of type
$\frac{1}{2}(1,1,1)$.

Let $T$ be the surface in $|-K_{X}|$, and let $S$ and $\grave{S}$
be general surfaces in $|-2K_{X}|$ that passes through the point
$P$. Then $S$ and $\grave{S}$ are irreducible and $S\supset
L\subset\grave{S}$, but $S\ne D\ne \grave{S}$.

Suppose now that $L$ is contained in $T$. Then $C=2L$ and
$-K_{X}\cdot L=2/5$, but the singularities of the curve $L$
consists of at most double points. Put $D\vert_{T}=mL+\Upsilon$,
where $m$~is~a~non-negative integer, and $\Upsilon$ is an
effective cycle whose support does not contain $L$. Then
$$
2n/5-4m/5=S\cdot\Upsilon\geqslant\mathrm{mult}_{P}\big(\Upsilon\big)\geqslant\mathrm{mult}_{P}\big(D\big)-\mathrm{mult}_{P}\big(L\big)>n-2m,
$$
which implies that $m>n/2$. But $m\leqslant n/2$ by
Remark~\ref{remark:curves}, which implies that $L\not\subset T$.

Suppose that $C=L$. Then $\mathrm{mult}_{P}(L)\leqslant 2$. Put
$D\cdot\grave{S}=\grave{m}C+\grave{\Upsilon}$, where
$\grave{m}$~is~a~non-negative integer and $\grave{\Upsilon}$ is an
effective cycle whose support does not contain $C$. Then
$$
4n/5-8\grave{m}/5=
S\cdot\grave{\Upsilon}\geqslant\mathrm{mult}_{P}\big(\grave{\Upsilon}\big)>n-2m,
$$
which implies that $m>n/2$. But $m\leqslant n/2$ by
Remark~\ref{remark:curves}, which implies that $C\ne L$.

The curve $L$ does not pass through a point of type
$\frac{1}{2}(1,1,1)$, and it follows from the generality of the
threefold $X$ that $C=L+Z$, where $Z$ is an irreducible curve such
that $Z\ne L$.

Put $D\big\vert_{S}=m_{L}L+m_{Z}Z+\Omega$, where $m_{L}$ and
$m_{Z}$ are non-negative integers and $\Omega$ is an effective
cycle whose support does not contain  $L$~and~$Z$.
We~may~assume~that $-K_{X}\cdot L\leqslant -K_{X}\cdot Z$, which
implies that either $-K_{X}\cdot L=1/5$ and $-K_{X}\cdot Z=3/5$,
or $-K_{X}\cdot L=-K_{X}\cdot Z=2/5$.

Suppose that $-K_{X}\cdot L=2/5$. Then $L$ and $Z$ are smooth
outside of $O$, and
$$
4n/5-4m_{L}/5-4m_{Z}/5=\grave{S}\big\vert_{S}\cdot\Omega\geqslant\mathrm{mult}_{P}\big(\Omega\big)>n-m_{L}-m_{C},
$$
which implies that $m_{L}+m_{C}>n$. But $m_{L}+m_{C}\leqslant n$
by Remark~\ref{remark:curves}.

Thus, we have $-K_{X}\cdot L=1/5$. The hypersurface $X$ can be
given by an equation
$$
w^{2}z+wg\big(x,y,z,t\big)+h\big(x,y,z,t\big)=0\subset\mathbb{P}\big(1,2,2,3,5\big)\cong\mathrm{Proj}\Big(\mathbb{C}[x,y,z,t,w]\Big),
$$
where $\mathrm{wt}(x)=1$, $\mathrm{wt}(y)=\mathrm{wt}(z)=2$,
$\mathrm{wt}(t)=3$, $\mathrm{wt}(w)=5$, and $g$ and $h$ are
quasihomogeneous polynomials  of degree $7$ and $12$,
respectively.

Let $R$ be a surface on the threefold $X$  that is cut out by
$z=0$, and let $\bar{R}$ and $\bar{L}$ be proper transforms of $R$
and $L$ on $U$, respectively. Then~$\bar{R}\cdot \bar{L}<0$, which
implies that $L\subset R\supset Z$,~and the curve $L$ is
contracted by the  projection $X\dasharrow\mathbb{P}(1,2,2,3)$ to
a point.

Let $\bar{Z}$ be the proper transform of the curve $Z$ on the
threefold $U$, let $\pi\colon\bar{R}\to R$ be a birational
morphism induced by $\alpha$, and let $\bar{E}$ be the curve on
the surface $\bar{R}$ that is contracted by $\pi$. Then
$$
\bar{L}^{2}=-1,\ \bar{L}\cdot\bar{Z}=\bar{L}\cdot\bar{E}=1,\ \bar{Z}^{2}=-1/3,\ \bar{E}^{2}=-35/6,\ \bar{Z}\cdot\bar{E}=4/3%
$$
on the surface $\bar{R}$, which implies that $L^{2}=-29/35$,
$L\cdot Z=43/35$, $Z^{2}=-1/35$ on the surface~$R$.

Suppose that $P\in L\cap Z$. Then $m_{L}+3m_{C}\leqslant 5n$ by
Remark~\ref{remark:curves}, but
$$
n/5+29m_{L}/35-43m_{Z}/35=\Omega\cdot L>n-m_{L}-m_{Z},\ 2n/5-43m_{L}/35+m_{Z}/35=\Omega\cdot Z>n-m_{L}-m_{Z},%
$$
which leads to a contradiction. Hence, either $L\ni P\not\in Z$,
or $Z\ni P\not\in L$.

Suppose that $Z\ni P\not\in L$. Then $\Omega\cdot Z>n-m_{Z}$ and
$\Omega\cdot L\geqslant 0$, which implies a contradiction.

Thus, we see that $L\ni P\not\in Z$. Then we have
$$
n/5+29m_{L}/35-43m_{Z}/35=\Omega\cdot L>n-m_{L},\
2n/5-43m_{L}/35+m_{Z}/35=\Omega\cdot Z\geqslant 0,
$$
which implies that $m_{L}<n$. Now it follows from Theorem~7.5 in
\cite{Ko97} that the log pair
$$
\left(R, L+\frac{m_{C}}{n}C+\frac{1}{n}\Omega\right)
$$
is not log
canonical at the point $P$. Then
$\mathrm{mult}_{P}(\Omega\vert_{L})>n$ by~Theorem~7.5 in
\cite{Ko97}, which implies that the inequality $\Omega\cdot L>n$
holds. The inequality $\Omega\cdot L>n$ leads to a contradiction.
\end{proof}

\begin{lemma}
\label{lemma:smooth-points-6} Suppose that $\gimel=6$. Then
$\mathrm{mult}_{P}(D)\leqslant n$.
\end{lemma}

\begin{proof}
Suppose that $\mathrm{mult}_{P}(D)>n$. It follows from
\cite{Ch06e} that the threefold $X$~has two quotient sin\-gu\-lar
points $O_{1}$ and $O_{2}$ of type $\frac{1}{2}(1,1,1)$ such that
there is a commutative~diagram
$$
\xymatrix{
&&&U\ar@{->}[dll]_{\sigma}\ar@{->}[d]^{\alpha}&&Y\ar@{->}[ll]_{\gamma}\ar@{->}[ddr]^{\eta}&&&\\%
&V\ar@{->}[dr]_{\omega}&&X\ar@{-->}[dl]_{\xi}\ar@{-->}[drrr]^{\psi}&&&&&\\
&&\mathbb{P}(1,1,1,4)\ar@{-->}[rrrr]_{\chi}&&&&\mathbb{P}(1,1,1),&&}
$$
where $\xi$, $\psi$ and $\chi$ are projections, $\alpha$ is a blow
up of $O_{1}$ with weights $(1,1,1)$, $\gamma$ is a  blow up with
weights $(1,1,1)$ of the point that dominates $O_{2}$, $\eta$ is
an elliptic fibration, $\omega$ is a double cover, and
$\sigma$~is~a~birational morphism that contracts $48$ irreducible
curves $C_{1},\ldots, C_{48}$.

The threefold $U$ contains $48$ curves $Z_{1},\ldots, Z_{48}$ such
that $\alpha(Z_{i})\cup\alpha(C_{i})$ is a  fiber of $\psi$ over
the~point $\psi(C_{i})$. Put $\bar{Z}_{i}=\alpha(Z_{i})$ and
$\bar{C}_{i}=\alpha(C_{i})$. Let $L$ be a fiber of the projection
$\psi$~that~passes through the point $P$, and let $T_{1}$ and
$T_{2}$ be general surfaces in $|-K_{X}|$ that contain $P$.

Suppose that $L$ is irreducible. Put $D\cdot T_{1}=mL+\Upsilon$,
where $m$ is non-negative integer and~$\Upsilon$~is an effective
cycle whose support does not contain $L$. Then $m\leqslant n$ by
Remark~\ref{remark:curves}. But
$$
n-m=D\cdot T_{1}\cdot T_{2}-mT_{2}\cdot L=T_{2}\cdot\Delta\geqslant\mathrm{mult}_{P}\big(\Delta\big)>n-m\mathrm{mult}_{P}(L),%
$$
which implies that $L$ is singular at the point $P$. Then there
is~a surface $T\in |-K_{X}|$ that~is~singular at the point $P$.
Let $S$ is a general surface in $|-2K_{X}|$ that contains $P$.
Then
$$
2n=D\cdot T\cdot S\geqslant\mathrm{mult}_{P}\big(D\cdot T\big)>2n,%
$$
which is a contradiction. Hence, the curve $L$ is reducible.

We have $L=\bar{C}_{i}\cup\bar{Z}_{i}$. Put
$D\vert_{T_{1}}=m_{1}\bar{C}_{i}+m_{2}\bar{Z}_{i}+\Delta$, where
$m_{1}$ and $m_{2}$ are non-negative integers and $\Delta$ is an
effective cycle whose support does not contain $\bar{C}_{i}$ and
$\bar{Z}_{i}$.

In the case when $P\in\bar{C}_{i}\cap\bar{Z}_{i}$, there is $T\in
|-K_{X}|$ such that $T$ singular at $P$, and we can obtain a
contradiction as above. So we may assume that $P\in\bar{C}_{i}$
and $P\not\in\bar{Z}_{i}$. Then
$$
n-m_{1}/2-m_{2}/2=\big(\bar{C}_{i}+\bar{Z}_{i}\big)\cdot\Delta\geqslant\mathrm{mult}_{P}\big(\Delta\big)>n-m_{1}\mathrm{mult}_{P}(\bar{C}_{i})=n-m_{1},%
$$
because $\bar{C}_{i}$ is smooth. Hence, we see that $m_{1}>m_{2}$.
But we have
$$
n-m_{1}\leqslant\Delta\cdot\bar{Z}_{i}=n/2-m_{1}\bar{C}_{i}\cdot\bar{Z}_{i}-m_{2}\bar{Z}_{i}^{2}=n/2-2m_{1}+3m_{2}/2<n/2-m_{1}/2,%
$$
which gives $m_{1}>m_{2}>n$. But $m_{1}+m_{2}\leqslant 2n$ by
Remark~\ref{remark:curves}. So we have a contradiction.
\end{proof}

Arguing as in the proofs of Lemmas~\ref{lemma:smooth-points-18}
and \ref{lemma:smooth-points-6}, we see that $(X,\frac{1}{n}D)$ is
log canonical at $P$~for
$$
\gimel\in\Big\{7,8,9,10,12,13,14,16,19,20,22,23,24,25,32,33,38\Big\},
$$
which completes the proof Lemma~\ref{lemma:smooth-points}.

\section{Singular points.}
\label{section:singular-points}

In this section we prove the assertion of
Lemma~\ref{lemma:singular-points}. Let us use the assumptions and
notation of Lemma~\ref{lemma:singular-points}. Suppose that $P$ is
a singular point of $U$. Let us derive a contradiction.

The point  $P$ is a singular point of type
$\frac{1}{\bar{r}}(1,\bar{a},\bar{r}-\bar{a})$, where $\bar{a}$
and $\bar{r}$ are coprime natural numbers such that
$\bar{r}>2\bar{a}$. Let $\beta\colon W\to U$ be a blow up~of $P$
with weights $(1,\bar{a},\bar{r}-\bar{a})$. Then
$$
-K_{W}^{3}=-K_{X}^{3}-\frac{1}{ra(r-a)}-\frac{1}{\bar{r}\bar{a}(\bar{r}-\bar{a})}.%
$$

Let $\breve{D}$ be the proper transform of $D$ on $W$. There is a
rational number $\nu$ such that
$$
\breve{D}\equiv (\alpha\circ\beta)^{*}\big(-nK_{X}\big)-\mu\beta^{*}\big(E\big)-\nu G,%
$$
where $G$ is the $\beta$-exceptional divisor. Then
$$
K_{W}+\frac{1}{n}\breve{D}+\big(\mu/n-1/r\big)\breve{E}\equiv
\beta^{*}\Big(K_{U}+\frac{1}{n}\bar{D}+\big(\mu/n-1/r\big)E\Big)-\epsilon G\equiv -\epsilon G,%
$$
where $\breve{E}$ is a proper transform of $E$ on the threefold
$W$, and $\epsilon\in\mathbb{Q}$. Then $\epsilon>0$ due to
\cite{Ka96}.

\begin{lemma}
\label{lemma:singular-points-K-cube-is-not-zero} The inequality
$-K_{W}^{3}\ne 0$ holds.
\end{lemma}

\begin{proof}
Suppose that $-K_{W}^{3}=0$. Then it follows from \cite{Ch06e}
that the linear system $|-rK_{W}|$ induces an elliptic fibration
$\eta\colon W\to Y$ for $r\gg 0$. Then $0\leqslant \breve{D}\cdot
C=-\epsilon G\cdot C<0$, where $C$ is a general fiber of the
elliptic fibration $\eta$.  So we have a contradiction.
\end{proof}

Thus, it follows from \cite{Ch06e} that either $-K_{W}^{3}<0$, or
$-K_{W}$ is nef and big.

\begin{lemma}
\label{lemma:singular-points-big-not-nef} Suppose that
$-K_{W}^{3}<0$. Then $-K_{W}$ is not big.
\end{lemma}

\begin{proof}
Suppose that $-K_{W}$ is big. Then it follows from \cite{Ch06e}
that we have the following possibilities:
\begin{itemize}
\item the equality $\gimel=25$ holds, and $O$ is a singular point of type $\frac{1}{7}(1,3,4)$;%
\item the equality $\gimel=43$ holds, and $O$~is a singular point of type $\frac{1}{9}(1,4,5)$;%
\end{itemize}
but both cases are similar. So we assume that $\gimel=43$. Then
$-K_{W}-4\beta^{*}(K_{U})$ is nef (see~\cite{Ch06e}), and there is
a surface $H$ in the linear system $|-2K_{X}|$ such that
$$
\breve{H}\equiv (\alpha\circ\beta)^{*}\big(-2K_{X}\big)-\frac{11}{9}\beta^{*}(E)-\frac{3}{2}G,\\%
$$
where $\breve{H}$ is a proper transform of the surface $H$ on the
threefold $W$. Thus, we have
$$
0\leqslant
\breve{H}\cdot\breve{D}\cdot\Big(-K_{W}-4\beta^{*}\big(K_{U}\big)\Big)=5n/9-11\mu/4-\nu,
$$
which is impossible, because $\nu-n/3+3\mu/4=n\epsilon>0$ and
$\mu>n/9$.
\end{proof}

Let $T$ be a surface in $|-K_{X}|$, and $\mathcal{P}$ be the
pencil generated by the divisors $nT$ and $D$. Then
\begin{equation}
\label{equation:crepant}
\mathcal{B}\equiv  -nK_{W}\equiv (\alpha\circ\beta)^{*}\big(-nK_{X}\big)-\frac{n}{r}\beta^{*}\big(E\big)-\frac{n}{\bar{r}} G,
\end{equation}
where $\mathcal{B}$ is the proper transforms of the pencil
$\mathcal{P}$ on the threefold $W$.

\begin{lemma}
\label{lemma:singular-points-nef-big} The divisor $-K_{W}$ is nef
and big.
\end{lemma}

\begin{proof}
Suppose that the divisor $-K_{W}$ is not nef and big. Then
$-K_{W}^{3}<0$, but $-K_{W}$ is not big by
Lemma~\ref{lemma:singular-points-big-not-nef}. Then the
equivalence~\ref{equation:crepant}~almost
uniquely~determines\footnote{For example, it follows from
\cite{ChPa05h} that the equivalence~\ref{equation:crepant} implies
that $n=1$ in the case when $a_{1}=1$.} the pencil $\mathcal{P}$
due to \cite{ChPa05h}.

All possible cases are similar. So we assume that
$\gimel\in\{45,48,58,69,74,79\}$. Then $O$~is~a~singular point of
type $\frac{1}{a_{4}}(1,a_{1},a_{3})$, and $X$ can be given by an
equation
$$
w^{2}z+wf\big(x,y,z,t\big)+g\big(x,y,z,t\big)=0\subset\mathbb{P}\big(1,a_{1},a_{2},a_{3},a_{4}\big)\cong\mathrm{Proj}\Big(\mathbb{C}[x,y,z,t,w]\Big),
$$
where $\mathrm{wt}(x)=1$, $\mathrm{wt}(y)=a_{1}$,
$\mathrm{wt}(z)=a_{2}$, $\mathrm{wt}(t)=a_{3}$,
$\mathrm{wt}(w)=a_{4}$, and $f$ and $g$~ are  polynomials.

Let $S$ be a surface that is cut out  on the threefold $X$  by
$z=0$, and $\mathcal{M}$ be a pencil generated by the divisors
$a_{2}T$ and $S$. Then it follows from \cite{ChPa05h} that either
$\mathcal{P}=\mathcal{M}$, or $\mathcal{P}=|-a_{1}K_{X}|$.

Suppose that $\mathcal{P}=|-a_{1}K_{X}|$. Then $\mu=n/a_{1}$,
which is impossible, because $\mu>n/a_{4}$.

We see that $\mathcal{P}=\mathcal{M}$. Let $M$ be a divisor in
$\mathcal{M}$, and $\bar{M}$ be its proper transform on $U$. Then
$$
\bar{M}\equiv \alpha^{*}\big(M\big)-\frac{a_{3}}{a_{4}}E
$$
in the case when $M\ne S$, but $\mu>n/a_{4}$. Thus, we see that
$D=S$, but $(X, \frac{1}{a_{2}} S)$ is log~ca\-no\-ni\-cal at the
point $O$ by Lemma~8.12 and Proposition~8.14~in~\cite{Ko97}, which
is a contradiction.
\end{proof}

Taking into account the possible values of
$(a_{1},a_{2},a_{3},a_{4})$, we see~that
$$
\gimel\in\Big\{8,12,13,16,20,24,25,26,31,33,36,38,46,47,48,54,56,58,65,74,79\Big\}
$$
by Lemmas~\ref{lemma:singular-points-K-cube-is-not-zero},
\ref{lemma:singular-points-big-not-nef} and
\ref{lemma:singular-points-nef-big} (see \cite{Ch06e}).

\begin{lemma}
\label{lemma:singular-points-20-25-31-33-38-58-nef} The case
$\gimel\not\in\{12,13,20,25,31,33,38,58\}$ is impossible.
\end{lemma}

\begin{proof}
Suppose that $\gimel\not\in\{12,13,20,25,31,33,38,58\}$. Then
$r=a_{4}$, $r-a=a_{3}$, $\bar{r}=r-a$, $\bar{a}=a$, and
$n\epsilon=\nu-(\bar{r}-\bar{a})(n/r-\mu)/\bar{r}-n/\bar{r}$. We
may assume that $\gimel\ne 24$, because the case $\gimel=24$~can
be considered in a similar way. Then $X$ can be given by the
equation
$$
w^{2}z+wf\big(x,y,z,t\big)+g\big(x,y,z,t\big)=0\subset\mathbb{P}\big(1,a_{1},a_{2},a_{3},a_{4}\big)\cong\mathrm{Proj}\Big(\mathbb{C}[x,y,z,t,w]\Big),
$$
where $\mathrm{wt}(x)=1$, $\mathrm{wt}(y)=a$,
$\mathrm{wt}(z)=d-2a_{4}$, $\mathrm{wt}(t)=a_{3}$,
$\mathrm{wt}(w)=a_{4}$, the point~$O$~is~given by the equations
$x=y=z=t=0$, and $f$ and $g$~are
quasi\-ho\-mo\-ge\-ne\-ous~polynomials. Then
$$
\breve{R}\equiv (\alpha\circ\beta)^{*}\big(-a_{2}K_{X}\big)-\frac{d-r}{r}\beta^{*}\big(E\big)-\frac{\bar{r}-\bar{a}}{\bar{r}} G,%
$$
where $\breve{R}$ is a proper transform on $W$ of the surface cut
out on  $X$ by $z=0$. Then $\breve{D}\ne\breve{R}$,~and
$$
n\frac{\sum_{i=1}^{4}a_{i}}{a_{1}a_{3}a_{4}}-\frac{\mu\big(d-r\big)}{a\big(r-a\big)}-\frac{\nu\big(\bar{r}-\bar{a}\big)}{\bar{a}\big(\bar{r}-\bar{a}\big)}=-K_{W}\cdot\breve{D}\cdot\breve{R}\geqslant 0,%
$$
which implies that $\mu<n/r$, because $\epsilon>0$. However we
know that $\mu>n/r$.
\end{proof}

So, the divisor $-K_{W}$ is nef and big, and
$\gimel\in\{12,13,20,25,31,33,38,58\}$, which implies that
$$
r=a_{4},\ r-a=a_{3},\ \bar{a}=a_{1},\ \bar{r}-\bar{a}=a_{2},\ a_{2}\ne a_{3},\ n\epsilon=\nu+(r-2a)\mu/(r-a)-2n/r%
$$
due to \cite{Ch06e}. Then $W$ has a singular point $\bar{P}\ne P$
of type $\frac{1}{\bar{r}}(1,\bar{a},\bar{r}-\bar{a})$ such that
the diagram
$$
\xymatrix{
U\ar@{->}[d]_{\alpha}&&W\ar@{->}[ll]_{\beta}&&V\ar@{->}[ll]_{\gamma}\ar@{->}[d]^{\eta}\\
X\ar@{-->}[rrrr]_{\psi}&&&&\mathbb{P}\big(1,a_{1},a_{2}\big)}
$$
commutes, where $\psi$ is a natural projection, $\gamma$ be a blow
up of the point $\bar{P}$ with weights
$(1,\bar{a},\bar{r}-\bar{a})$, and $\eta$ is an elliptic
fibration. Let $F$ be the exceptional divisor of $\gamma$, and
$\bar{G}$ be the proper transform of the divisor $G$ on the
threefold $V$. Then $F$ and $\bar{G}$ are sections of $\eta$, and
$G\not\ni\bar{P}\not\in \breve{E}$.

It follows from the inequality $-K_{W}\cdot \breve{D}\geqslant 0$
and the proof of Lemma~\ref{lemma:discrepancy} that $\epsilon<1$,
which implies that the log pair $(W,
\frac{1}{n}\breve{D}+(\mu/n-1/r)\breve{E}+\epsilon G)$ is not log
canonical at some point $Q\in G$.

\begin{lemma}
\label{lemma:singular-points-singular-point} The threefold $W$ is
smooth at the point $Q$.
\end{lemma}

\begin{proof}
Suppose that $W$ is singular at the point $Q$. Then $Q$ is a
singular point of type $\frac{1}{\breve{r}}(1,1,\breve{r}-1)$,
where either $\breve{r}=\bar{r}-\bar{a}$, or $\breve{r}=\bar{a}\ne
1$. Let $\omega\colon \breve{W}\to W$ be a blow up of $Q$ with
weights~$(1,1,\breve{r}-1)$, and $\mathcal{H}$ be the proper
transform of $\mathcal{P}$ on $\breve{W}$. Then it follows from
\cite{Ka96} that $\mathcal{H}\equiv -nK_{\breve{W}}$, which
implies that $n=r\mu=a_{1}$ due  to \cite{ChPa05h}. However we
know that $\mu>n/r$.
\end{proof}

Thus, it follows from
Lemma~\ref{lemma:singular-points-singular-point} that
$\mathrm{mult}_{Q}(\breve{D})>n-n\epsilon-(\mu-n/r)\mathrm{mult}_{Q}(\breve{E})$.

\begin{lemma}
\label{lemma:singular-points-T} There is a surface $T\in |-K_{W}|$
such that $Q\in T$.
\end{lemma}

\begin{proof}
The existence of a surface $T\in |-K_{W}|$ that passes through the
point $Q$ is obvious in the~case when $a_{1}=1$. Thus, we may
assume that $a_{1}\ne 1$. Then $\gimel\in\{33,38,58\}$, but we
consider only the case $\gimel=38$, because the cases $\gimel=33$
and $\gimel=58$ can be considered in a similar way.

Suppose that $\gimel=38$. Then there is a unique surface $T\in
|-K_{W}|$. Suppose that $Q\not\in T$.

Arguing as in the proof of Lemma~\ref{lemma:lower-bound-for-mu},
we see that $\mathrm{mult}_{Q}(\breve{D})\leqslant
(a_{1}+a_{2})\nu/a_{1}$. Then
$$
\nu\frac{a_{1}+a_{2}}{a_{1}}>n-\big(\mu-n/7\big)-\big(\nu+3\mu/5-2n/7\big),
$$
but $\mathrm{mult}_{Q}(\breve{D})>n+n/r-\mu-n\epsilon$. Thus, we
have $\mu>55n/56-5\nu/2$.

The inequality $-K_{W}\cdot \breve{D}\geqslant 0$ and the proof of
Lemma~\ref{lemma:discrepancy} imply that $\nu\leqslant 10\mu/7$
and $\mu\leqslant 9n/40$, respectively. The hypersurface $X$ can
be given by the equation
$$
w^{2}y+w\Big(t^{2}+tf_{5}\big(x,y,z\big)+f_{10}\big(x,y,z\big)\Big)+tf_{13}\big(x,y,z\big)+f_{18}\big(x,y,z\big)=0\subset\mathrm{Proj}\Big(\mathbb{C}[x,y,z,t,w]\Big),
$$
where $\mathrm{wt}(x)=1$, $\mathrm{wt}(y)=2$, $\mathrm{wt}(z)=3$,
$\mathrm{wt}(t)=5$, $\mathrm{wt}(w)=8$, and $f_{i}(x,y,z)$ is a
quasi\-ho\-mo\-ge\-ne\-ous polynomial of degree $i$. Let
$\breve{S}$ be the proper transform on the threefold $W$ of the
surface that is cut out on $X$ by the equation
$wy+(t^{2}+tf_{5}(x,y,z)+f_{10}(x,y,z))=0$. Then
$$
\breve{S}\equiv \big(\alpha\circ\beta\big)^{*}\big(-10K_{X}\big)-\frac{18}{8}\beta^{*}\big(E\big)-\frac{13}{5}G,%
$$
but $\breve{S}\ne\breve{D}$. The divisor $-K_{W}$ is nef. Hence,
we have
$$
0\leqslant -K_{W}\cdot\breve{D}\cdot\breve{S}=3n/4-6\mu/5-13\nu/6,
$$
but $\nu\leqslant 8\mu/5$, which implies $\nu\leqslant 9n/35$. Now
we can easily obtain a contradiction.
\end{proof}

It follows from \cite{Ch06e} that $|-rK_{W}|$ does not have base
points for $r\gg 0$ and induces a birational morphism
$\omega\colon W\to\bar{W}$ such that $\bar{W}$ is a hypersurface
in $\mathbb{P}(1,a_{1},a_{2},2a_{3},3a_{3})$ of degree $6a_{3}$

\begin{lemma}
\label{lemma:singular-points-omega} The morphism $\omega$ is not
an isomorphism in a neighborhood of the point $Q$.
\end{lemma}

\begin{proof}
Suppose that $\omega$ is an isomorphism in a neighborhood of the
point $Q$. Then it follows from the proof of Theorem~5.6.2 in
\cite{CPR} that there is $R\in|-s2a_{1}a_{3}K_{W}|$ such that
$\mathrm{mult}_{Q}(R)\geqslant s$, but $R$ does~not contain
components of the cycle $\breve{D}\cdot S$ that pass through $Q$,
where $s\in\mathbb{N}$. Then
$$
s2a_{3}\big(-na_{1}K_{X}^{3}-\mu/a_{3}-\nu/a_{2}\big)=R\cdot\breve{D}\cdot
T\geqslant\mathrm{mult}_{Q}\big(\breve{D}\cdot T\big)s
>\big(n-\nu-\mu(a_{3}-a_{1})/a_{3}+2n/a_{4}\big)s,%
$$
because $Q\not\in\breve{E}$. Now we can derive a contradiction
using $n\epsilon=\nu+(a_{3}-a_{1})\mu/a_{3}-2n/a_{4}>0$.
\end{proof}

It follows from Lemma~\ref{lemma:singular-points-omega} that there
is a unique curve $C\subset W$~that contains $Q$ such~that
$$
-K_{W}\cdot C=0,\ \beta^{*}\big(-K_{U}\big)\cdot C=1/a_{4},\ C\cdot G=1,%
$$
which implies that $\gimel\not\in\{33,38,58\}$ by
Lemma~\ref{lemma:singular-points-T}. Hence, we have
$\gimel\in\{12,13,20,25,31\}$.

Put $\breve{D}\cdot T=mC+\Omega$, where $m$ is a non-negative
integer, and $\Omega$ is an effective one-cycle, whose support
does not contain the curve $C$. Then it follows from
Remark~\ref{remark:curves} that
$$
m\leqslant 5n/4-\mu,\ m\leqslant 11n/15-\mu/2,\ m\leqslant 13n/15-\mu,\ m\leqslant 5n/7-\mu/3,\ m\leqslant 2n/3-\mu%
$$
in the case when $\gimel=12,13,20,25,31$, respectively. Recall
that $\bar{G}$ is a section of $\eta$.

Let $\mathcal{H}$ be a pencil in $|-a_{2}K_{W}|$ of surfaces
passing through the point $Q$, and $H$ be a general surface in
$\mathcal{H}$. Then $C$ is the only curve in the base locus of
$\mathcal{H}$ that passes through $Q$. Then
$$
-na_{2}K_{X}^{3}-a_{2}\mu/(a_{1}a_{3})-\nu/a_{1}=H\cdot\Omega\geqslant\mathrm{mult}_{Q}\big(\Omega\big)>n-\nu-\mu(a_{3}-a_{1})/a_{3}+2n/a_{4}-m,
$$
which immediately implies that either $\gimel=12$, or $\gimel=13$.

\begin{lemma}
\label{lemma:singular-points-12} The inequality $\gimel\ne 12$
holds.
\end{lemma}

\begin{proof}
Suppose that $\gimel=12$. Let $R$ be a sufficiently general
surface in $|-2K_{W}|$ that contains the~point $Q$. Then
$R\vert_{T}=C+L+Z$, where $L=G\vert_{T}$, the curve $Z$ is
reduced, and $P\not\in\beta(Z)$.

Suppose that $Z$ is irreducible. Then $Z^{2}=-4/3$, $C^{2}=-2$,
$L^{2}=-3/2$ on the surface $T$. Put
$$
\breve{D}\vert_{T}=m_{C}C+m_{L}L+m_{Z}Z+\Upsilon,
$$
where $m_{C}$, $m_{L}$ and $m_{Z}$ are non-negative integers, and
$\Upsilon$ is an effective one-cycle, whose support does not
contain the curve $C$, $L$ and $Z$.

Suppose that $Q\not\in\breve{E}$. Then $m_{C}>2n/3-m_{Z}/3$,
because
$$
5n/6-2\mu/3-\nu=R\cdot\breve{D}\cdot T=m_{L}+m_{Z}/3+R\cdot\Upsilon>m_{L}+m_{Z}/3+3n/2-\nu-2\mu/3-m_{L}-m_{C},%
$$
but $4m_{Z}/3\geqslant 2m_{C}-n/3$, because $\Upsilon\cdot
Z\geqslant 0$. Thus, we have
$$
m_{C}>2n/3+m_{Z}/3\geqslant 7n/12+m_{C}/2,
$$
which gives $m_{C}>7n/6$, but $m_{C}\leqslant 5n/6$ by
Remark~\ref{remark:curves}, because
$-K_{X}\cdot\alpha\circ\beta(C)=5/6$.

Thus, we see that $Q\in\breve{E}$. Then $C\subset\breve{E}$ and
$\beta(C)\in|\mathcal{O}_{\mathbb{P}(1,\,1,\,3)}(1)|$, but
$$
5n/6-2\mu/3-\nu=R\cdot\breve{D}\cdot T=m_{L}+m_{Z}/3+R\cdot\Upsilon>m_{L}+m_{Z}/3+7n/4-\nu-5\mu/3-m_{L}-m_{C},%
$$
which gives $m_{C}>11n/12-\mu+m_{Z}/3$. We have
$-K_{X}\cdot\alpha\circ\beta(Z)=5/6$ and $Z\cdot\breve{E}=2$, but
$$
4m_{Z}/3\geqslant 2m_{C}+2\mu-5n/6,
$$
because $Z\cdot\Upsilon\geqslant 0$. Thus, we have $m_{Z}>3n/2$,
but $m_{Z}\leqslant n/2$ by Remark~\ref{remark:curves}.

Therefore, the curve $Z$ is reducible. Then $Q\in\breve{E}$ and
$Z=\acute{Z}+\grave{Z}$, where $\acute{Z}$ and
$\grave{Z}$~are~irreducible curves such that
$G\cdot\acute{Z}=G\cdot\grave{Z}=-K_{U}\cdot\beta(\grave{Z})=0$
and $-K_{X}\cdot\alpha\circ\beta(\acute{Z})=7/12$. Then
$$
\acute{Z}^{2}=-4/3,\ \grave{Z}^{2}=C^{2}=-2,\ L^{2}=-3/2,\ L\cdot C=\acute{Z}\cdot C=\acute{Z}\cdot\grave{Z}=\grave{Z}\cdot C=1,\ L\cdot\grave{Z}=L\cdot\acute{Z}=0%
$$
on the surface $T$. Put
$\breve{D}\vert_{T}=\bar{m}_{C}C+\bar{m}_{L}L+\bar{m}_{Z}\acute{Z}+\Phi$,
where $\bar{m}_{C}$, $\bar{m}_{L}$, $\bar{m}_{Z}$ are non-negative
integers, and $\Phi$ is an effective cycle, whose support does not
contain $C$, $L$ and $\acute{Z}$. Then
$$
R\big\vert_{T}\cdot\Phi\geqslant\mathrm{mult}_{Q}\big(\Phi\big)>7n/4-\nu-5\mu/3-m_{L}-m_{C},
$$
and $\Phi\cdot \acute{Z}\geqslant 0$. We have
$\beta^{*}(-K_{U})\vert_{T}\cdot\Phi\geqslant 0$. Thus, we see
that
$$
\bar{m}_{C}>11n/12-\mu+\bar{m}_{Z}/3,\ 4\bar{m}_{Z}/3\geqslant \bar{m}_{C}+\mu-5n/6,\ \bar{m}_{C}+\mu\leqslant 5/4-\bar{m}_{Z},%
$$
but this system of linear inequalities is inconsistent, which
completes the proof.
\end{proof}

Thus, we see that $\gimel=13$. Then $C\subset\breve{E}$, because
otherwise we have
$$
2\Big(11n/30-\mu/6-\nu/2\Big)=H\cdot\Omega\geqslant\mathrm{mult}_{Q}\big(\Omega\big)>7n/5-\nu-\mu/3-m,%
$$
which implies that $m>2n/3$, which is impossible, because
$m\leqslant 11n/15-\mu/2$ and $\mu>n/5$.~Put
$$
\bar{D}\big\vert_{\breve{E}}=\bar{m}C+\Upsilon,
$$
where $\bar{m}$ is a non-negative integer, and $\Upsilon$ is an
effective cycle, whose support does not contain the curve $C$.
Then $\bar{m}\leqslant 5\mu/2$, because we have
$\beta(C)\in|\mathcal{O}_{\mathbb{P}(1,\,2,\,3)}(2)|$ and the
curve $C$ is reduced, where $E\cong\mathbb{P}(1,2,3)$. Then
$\bar{m}\leqslant 11n/12$, because $\mu\leqslant 11n/30$.

The log pair $(W,\ \frac{1}{n}\breve{D}+\breve{E}+\epsilon G)$ is
not log canonical at the point $Q$. Hence, the log pair
$$
\left(\breve{E},\ C+\frac{\nu+\mu/3-2n/5}{n}G\big\vert_{\breve{E}}+\frac{1}{n}\Upsilon\right)%
$$
is not log canonical at the point $Q$ by Theorem~7.5 in
\cite{Ko97}. Then
$$
5\mu/3-\nu=\big(\bar{m}C+\Upsilon\big)\cdot C=\Upsilon\cdot C>7n/5-\nu-\mu/3,%
$$
because $\mathrm{mult}_{Q}(\Upsilon\vert_{C})>7n/5-\nu-\mu/3$ by
Theorem~7.5 in \cite{Ko97}. Thus, we see that $\mu>7n/10$,~which
is impossible, because  $\mu\leqslant 11n/30<7n/10$. The assertion
of Lemma~\ref{lemma:singular-points} is proved.

\section{Quadratic involutions.}
\label{section:quadratics-involutions}

In this section we prove the assertion of
Lemma~\ref{lemma:quadratic-involutions}. Let us use the
assumptions and notation of
Lemma~\ref{lemma:quadratic-involutions}. Suppose that
$d=2r+a_{j}$. To prove Lemma~\ref{lemma:quadratic-involutions} we
must derive a contradiction.

It follows from the equality $d=2r+a_{j}$ that the threefold $X$
can be given by the equation
$$
x_{i}^{2}x_{j}+x_{i}f\big(x_{0},x_{1},x_{2},x_{3},x_{4}\big)+g\big(x_{0},x_{1},x_{2},x_{3},x_{4}\big)=0\subset\mathrm{Proj}\Big(\mathbb{C}[x_{0},x_{1},x_{2},x_{3},x_{4}]\Big),
$$
where $i\ne j$, $a_{i}=r$, $\mathrm{wt}(x_{0})=1$,
$\mathrm{wt}(x_{k})=a_{k}$, $f$ and $g$ are
quasi\-ho\-mo\-ge\-ne\-ous polynomials that do not depend on
$x_{i}$. Put $\bar{a}_{3}=a_{3+4-i}$, $\bar{a}_{4}=a_{i}a_{j}$,
$\bar{d}=2\bar{a}_{4}$. Then there is a commutative~diagram
$$
\xymatrix{
&U\ar@{->}[d]_{\sigma}\ar@{->}[rrrr]^{\alpha}&&&&X\ar@{-->}[d]^{\xi}\\
&V\ar@{^{(}->}[rr]&&\mathbb{P}\big(1,a_{1},a_{2},\bar{a}_{3},\bar{a}_{4}\big)\ar@{-->}[rr]_{\chi}&&\mathbb{P}\big(1,a_{1},a_{2},\bar{a}_{3}\big),}
$$
where $\xi$ and $\chi$ are projections, and $\sigma$ is a
birational morphism that contracts smooth irreducible rational
curves $C_{1},\ldots, C_{l}$ such that $V$ is a hypersurface in
$\mathbb{P}(1,a_{1},a_{2},\bar{a}_{3},\bar{a}_{4})$ of degree
$\bar{d}$~with terminal non-$\mathbb{Q}$-factorial singularities,
where $l=a_{i}a_{j}(d-a_{i})\sum_{i=1}^{4}a_{i}$. Then
$-K_{X}\cdot\alpha(C_{k})=1/a_{i}$.

Let $M$ be the surface that is cut out on the threefold $X$ by
$x_{j}=0$, and $\bar{M}$ be the proper transform of $M$ on the
threefold $U$. Then $M\ne D$ by Lemma~8.12 and Proposition~8.14 in
\cite{Ko97}.

\begin{lemma}
\label{lemma:inequality-for-mu} The inequalities $\mu\leqslant
-a_{j}nK_{X}^{3}(r-a)a/(d-r)\leqslant n(d-r)/(ra_{j})$ hold.
\end{lemma}

\begin{proof}
The divisor $-K_{U}$ is nef. The inequality $\mu\leqslant
-a_{j}nK_{X}^{3}(r-a)a/(d-r)$ follows from
$$
0\leqslant -K_{U}\cdot\bar{M}\cdot\bar{D}=-a_{j}nK_{X}^{3}-\mu\big(d-r\big)/\big(ar-a^{2}\big),%
$$
and to conclude the proof we must show that
$-a_{j}nK_{X}^{3}(r-a)a/(d-r)\leqslant n(d-r)/(ra_{j})$.
\newpage
Suppose that $-a_{j}nK_{X}^{3}(r-a)a/(d-r)>n(d-r)/(ra_{j})$. Then
$$
\frac{d-r}{ra_{j}}<-a_{j}K_{X}^{3}\frac{\big(r-a\big)a}{d-r}=\frac{da_{j}(r-a)a}{\big(d-r\big)a_{1}a_{2}a_{3}a_{4}},
$$
but $a_{1}a_{2}a_{3}a_{4}\geqslant a_{j}r(r-a)a$. Thus, we have
$(d-r)^{2}<d(d-2r)$, which is a contradiction.
\end{proof}

We have $E\cong\mathbb{P}(1,a,r-a)$, and
$|\mathcal{O}_{\mathbb{P}(1,\,a,\,r-a)}(1)|$ consists of a single
curve when $a\ne 1$.

\begin{lemma}
\label{lemma:a-not-1} The inequality $a\ne 1$ holds.
\end{lemma}

\begin{proof}
Suppose that $a=1$. Taking into account the possible values of
$(a_{1},a_{2},a_{3},a_{4})$, we see~that
$$
\gimel\in\Big\{6,7,8,9,12,13,16,15,17,20,25,26,30,36,31,41,47,54\Big\},
$$
but we only consider the cases $\gimel=7$ and $\gimel=36$. The
remaining $16$ cases can be~considered in a similar way. So the
reader can easily obtain a contradiction in these cases by
himself.

Suppose that $\gimel=7$. Then $X$ is a hypersurface in
$\mathbb{P}(1,1,2,2,3)$ of degree $8$, which implies~that the
point $O$ is a singular point of type $\frac{1}{3}(1,1,2)$. Let
$S$ be the unique surface in $|-K_{U}|$ that contains the point
$P$. Then $S$ is an irreducible surface, which is smooth at the
point $P$.

The singularities of $U$ consists of singular points $P_{0}$,
$P_{1}$, $P_{2}$, $P_{3}$ and $P_{4}$ of type $\frac{1}{2}(1,1,1)$
such that $P_{0}$ is a singular point of $E$. It follows from
\cite{Ch06e} that there is a commutative diagram
$$
\xymatrix{
&U\ar@{->}[d]_{\alpha}&&Y_{i}\ar@{->}[ll]_{\beta_{i}}\ar@{->}[drr]^{\eta_{i}}&&&\\
&X\ar@{-->}[rrrr]_{\xi_{i}}&&&&\mathbb{P}(1,1,2),&}
$$
where $\xi_{i}$ is a projection, $\beta_{i}$ is a blow of $P_{i}$
with weights $(1,1,1)$, and $\eta_{i}$ is a morphism.

Suppose that $P\not\in \cup_{i=1}^{l}C_{i}$. The proper transform
of $E$ on the threefold $Y_{i}$ is a section of $\eta_{i}$ in
the~case when $i\ne 0$. Hence, there is a~surface~$H\in |-2K_{U}|$
such that
$$
2\big(2n/3-\mu/2\big)=\bar{D}\cdot H\cdot S\geqslant \mathrm{mult}_{P}\big(\bar{D}\big)>4n/3-\mu,%
$$
which is a contradiction. So we may assume that $P\in C_{1}$. Then
$-K_{X}\cdot\alpha(C_{1})=1/3$.

Let $Z_{1}$ be an irreducible curve such that
$\bar{M}\vert_{S}=C_{1}+Z_{1}$. Put $L=E\vert_{S}$. Then
$$
C_{1}^{2}=-2,\ Z_{1}^{2}=L^{2}=-3/2,\ C_{1}\cdot Z_{1}=L\cdot C_{1}=1,\ L\cdot Z_{1}=3/2%
$$
on the surface $S$. Put
$\bar{D}\vert_{S}=m_{C}C_{1}+m_{Z}Z_{1}+m_{L}L+\Omega$, where
$m_{C}$, $m_{Z}$ and $m_{L}$ are non-negative integers, and
$\Omega$ is an effective cycle, whose support does not contain
$C_{1}$, $Z_{1}$ and $L$. Then
$$
n-3\mu/2+3m_{Z}/2-3m_{L}/2-m_{C}=Z_{1}\cdot\Omega\geqslant 0,\ 3\mu/2-3m_{Z}/2+3m_{L}/2-m_{C}=L\cdot\Omega\geqslant 0,%
$$
which implies that $3m_{Z}/2\geqslant 3(\mu+m_{L})/2+m_{C}-n$ and
$3(\mu+m_{L})/2\geqslant 3m_{Z}/2+m_{C}$. Then
$$
4n/3-\mu-m_{L}-m_{Z}=\big(L+C_{1}+Z_{1}\big)\cdot\Omega\geqslant\mathrm{mult}_{P}\big(\Omega\big)>4n/3-\mu-m_{L}-m_{C},
$$
which gives $m_{C}>m_{Z}$ and $4n/3\geqslant \mu+m_{L}+m_{Z}$. So
we see that $m_{Z}\leqslant n/2$ and $m_{C}\leqslant n/2$.

Then it follows from Theorem~7.5 in \cite{Ko97} that the log pair
$$
\left(S,\ C_{1}+\frac{\bar{m}_{L}+\mu-n/3}{n}L+\frac{m_{Z}}{n}Z+\frac{1}{n}\Omega\right)%
$$
is not log canonical at $P$, because $m_{C}\leqslant n$. So it
follows from Theorem~7.5 in \cite{Ko97} that
$$
C_{1}\cdot\Omega\geqslant\mathrm{mult}_{P}\Big(\Omega\big\vert_{C_{1}}\Big)>n-m_{L}-\mu+n/3,
$$
which implies that $m_{C}>m_{Z}/2+n/2$, but $m_{C}\leqslant n/2$.
So the case $\gimel=7$ is impossible.

Now we suppose that $\gimel=36$. Then $X$ is a general
hypersurface in $\mathbb{P}(1,1,4,6,7)$ of degree $18$,
and~$O$~is~a~singular point of type $\frac{1}{7}(1,1,6)$. Arguing
as in the case $\gimel=7$, we see that $P\not\in
\cup_{i=1}^{l}C_{i}$, which implies that we may assume that $P\in
C_{1}$. Put $L=C_{1}$.

Let $S$ be a surface in $|-K_{U}|$ such that $P\in S$. Then
$\bar{M}\vert_{S}=L+Z$, where $Z$ is an irreducible curve. Put
$C=E\vert_{S}$. Then the intersection form of $C$, $L$, $Z$ on $S$
is given~by
$$
Z^{2}=C^{2}=-7/6,\ L^{2}=-2,\ Z\cdot L=C\cdot L=1,\ Z\cdot C=5/6,
$$
and $P$ is the intersection point of the curves $L$ and $C$. Put
$$
\bar{D}\big\vert_{S}=m_{L}L+m_{C}C+m_{Z}Z+\Omega,
$$
where $m_{L}$, $m_{C}$ and $m_{Z}$ are a non-negative integers,
and $\Omega$ is an effective cycle, whose support does not contain
the curves $L$, $C$, and $Z$.

It follows from the proof of Theorem~5.6.2 in \cite{CPR} that we
can find $H\in|-s6K_{U}|$ that has multiplicity at least $s>0$ at
the point $P$, but does~not contain~com\-po\-nents of $\Omega$
that pass through the point $P$, where $s$ is a natural number.
Then
$$
s6\big(3n/28-\mu/6-m_{C}/6-m_{Z}/6\big)=H\big\vert_{S}\cdot \Omega\geqslant \mathrm{mult}_{P}\big(\Omega\big)s>\big(8n/7-\mu-m_{L}-m_{C}\big) s,%
$$
which implies that $m_{L}>n/2+m_{Z}$, but $m_{L}\leqslant 3n/4$ by
Remark~\ref{remark:curves}. We have
$$
3n/28-\mu/6=-K_{U}\big\vert_{S}\cdot\Big(m_{L}L+m_{C}C+m_{Z}Z+\Omega\Big)\geqslant
-K_{U}\big\vert_{S}\cdot\Big(m_{L}L+m_{C}C+m_{Z}Z\Big)=\frac{m_{C}+m_{Z}}{6},
$$
which implies that $m_{C}+m_{Z}\leqslant 9n/14-\mu$. On the
surface $S$ we have
$$
7\mu/6+7m_{C}/6-5m_{Z}/6-m_{L}=\Omega\cdot C>8n/7-\mu-m_{L}-m_{C},%
$$
which implies that $13(\mu+m_{C})/6>8n/7+5m_{Z}/6$. The inequality
$\Omega\cdot Z\geqslant 0$ implies that
$$
2n/7-5\mu/6-m_{L}-5m_{C}/6+7m_{Z}/6\geqslant 0,
$$
which implies that $7m_{Z}/6\geqslant 5\mu/6+m_{L}+5m_{C}/6-2n/7$,
but $m_{Z}\leqslant 3n/8$ by Remark~\ref{remark:curves}.

It follows from Lemma~\ref{lemma:inequality-for-mu} that
$18n/77\geqslant\mu>n/7$. The  inequalities obtained
$$\left\{\aligned
&13\big(\mu+m_{C}\big)/6>8n/7+5m_{Z}/6,\\
&21n/48\geqslant 7m_{Z}/6\geqslant 5\mu/6+m_{L}+5m_{C}/6-2n/7,\\
&m_{C}+m_{Z}\leqslant 9n/14-\mu,\\
&3n/4\geqslant m_{L}>n/2+m_{Z},\\
&18n/77\geqslant\mu>n/7,\\
\endaligned
\right.
$$
are inconsistent.  So we have a contradiction. Thus, the case
$\gimel=36$ is impossible as well.
\end{proof}

Taking into account the possible values of the quadruple
$(a_{1},a_{2},a_{3},a_{4})$, we see that
$$
\gimel\in\big\{13,18,23,24,27,32,38,40,42,43,44,45,46,48,56,58,60,61,65,68,69,74,76,79\big\}
$$
by Lemmas~\ref{lemma:a-not-1}. Let $T$ be a general surface in
$|-K_{U}|$. Then $T\vert_{E}\in
|\mathcal{O}_{\mathbb{P}(1,\,a,\,r-a)}(1)|$.

\begin{lemma}
\label{lemma:smooth-points-upstairs-in-T} The point $P$ is
contained in the surface $T$.
\end{lemma}

\begin{proof}
It follows from Lemmas~\ref{lemma:lower-bound-for-mu} and
\ref{lemma:inequality-for-mu} that $P\in T$ unless
$\gimel\in\{13,24\}$. Therefore, we may assume that
$\gimel\in\{13,24\}$ and $P\not\in T$. Let us derive a
contradiction.

Let $L$ be the curve in
$|\mathcal{O}_{\mathbb{P}(1,\,a,\,r-a)}(1)|$. Then $P\not\in L$,
because $P\not\in T$. Thus, there is a unique smooth irreducible
curve $C$ in the linear system
$|\mathcal{O}_{\mathbb{P}(1,\,a,\,r-a)}(a)|$ that contains the
point $P$. Put
$$
\bar{D}\big\vert_{E}=\delta C+\Upsilon\equiv r\mu L,
$$
where $\delta$ is a non-negative integer, and $\Upsilon$ is an
effective cycle such that $C\not\subset\mathrm{Supp}(\Upsilon)$.

Arguing as in the proof of  Lemma~\ref{lemma:lower-bound-for-mu},
we see that $\delta\leqslant r\mu/a$, which gives $\delta<n$ by
Lemma~\ref{lemma:inequality-for-mu}.

It follows from Theorem~7.5 in \cite{Ko97} that $(E,
\frac{1}{n}\bar{D}\vert_{E})$ is not log canonical at $P$, which
implies that the log pair $(E, C+\frac{1}{n}\Upsilon)$ is not log
canonical at $P$. It follows from Theorem~7.5 in \cite{Ko97}~that
$$
r\mu/\big(r-a\big)\geqslant\big(r\mu-a\delta\big)/\big(r-a\big)=C\cdot\Upsilon\geqslant\mathrm{mult}_{P}\big(\Upsilon\vert_{C}\big)>n,%
$$
which implies that $\mu\geqslant n(r-a)/r$, which is impossible by
Lemma~\ref{lemma:inequality-for-mu}.
\end{proof}

It follows from \cite{CPR} that $T\cap E\cap
\cup_{i=1}^{l}C_{i}\ne\varnothing\iff\gimel\in\{43, 46, 69, 74,
76, 79\}$.

\begin{lemma}
\label{lemma:smooth-points-13-24-32-43-46} The case
$\gimel\not\in\{13, 24, 32, 43, 46\}$ is impossible.
\end{lemma}

\begin{proof}
Suppose that $\gimel\not\in\{13, 24, 32, 43, 46, 56\}$. It follows
from the proof of Theorem~5.6.2 in \cite{CPR} that~there are
$s\in\mathbb{N}$ and $H\in |-sa_{1}\bar{a}_{3}K_{U}|$ such that
$\mathrm{mult}_{P}(H)\geqslant s$, but $H$ does not contain
components of the cycle $\bar{D}\cdot T$ passing through $P$ that
are different from the curves $C_{1},\ldots,C_{l}$.

We have $\gimel\in\{69, 74, 76, 79\}$ and
$P\in\cup_{i=1}^{l}C_{i}$, because otherwise we get a
contradiction~using
$$
sa_{1}\bar{a}_{3}\left(-nK_{X}^{3}-\frac{\mu}{a\big(r-a\big)}\right)=\bar{D}\cdot H\cdot T\geqslant \mathrm{mult}_{P}\big(\bar{D}\big)s>\big(n+n/r-\mu\big) s.%
$$

We may assume that $P\in C_{1}$.~Put $\bar{D}\cdot
T=mC_{1}+\Delta$, where $m$ is a non-negative integer~number, and
$\Delta$ is an effective cycle, whose support does not contain the
curve $C_{1}$. Then it follows from the proof of Theorem~5.6.2 in
\cite{CPR} that there is a~surface $R\sim
-sa_{1}\bar{a}_{3}K_{U}$~such~that
$$
sa_{1}\bar{a}_{3}\left(-nK_{X}^{3}-\frac{\mu}{a\big(r-a\big)}\right)=R\cdot\Delta\geqslant\mathrm{mult}_{P}\big(\Delta\big)s>\big(n+n/r-\mu-m\big)s,%
$$
where $s\in\mathbb{N}$. The inequality obtained is impossible,
because $m\leqslant -a_{i}nK_{X}^{3}$ by
Remark~\ref{remark:curves}.

Suppose that $\gimel=56$. As in the previous case, there is $H\in
|-s24K_{U}|$ such that
$$
s24\big(n/22-\mu/24\big)=\bar{D}\cdot H\cdot T\geqslant \mathrm{mult}_{P}\big(\bar{D}\big)s>\big(12n/11-\mu\big) s,%
$$
where $s$ is a natural number. Now we can easily obtain a
contradiction with $\mu>n/r$.
\end{proof}

Thus, to complete the proof of
Lemma~\ref{lemma:quadratic-involutions}, we have to consider the
cases $\gimel=13, 24, 32, 43, 46$ one by one. For the sake of
simplicity, we only consider the cases $\gimel=13$ and
$\gimel=43$, because the remaining cases can be considered in a
similar way.

\begin{lemma}
\label{lemma:43} The inequality $\gimel\ne 43$ holds.
\end{lemma}

\begin{proof}
Suppose that $\gimel=43$. Then $X$ is a general hypersurface in
$\mathbb{P}(1,2,3,5,9)$ of degree $20$, and~$O$~is a singular
point of type $\frac{1}{9}(1,4,5)$. The base locus of $|-2K_{U}|$
consists of two irreducible curves $C$ and $L$ such that $L=T\cdot
E$, and $C$ is the curve among $C_{1},\ldots,C_{l}$ such that
$C\cap L\ne\varnothing$.

Suppose that $P\not\in C$. Then it follows from the proof of
Theorem~5.6.2 in \cite{CPR} that we can find a~surface~$H\in
|-s20K_{U}|$ that has multiplicity at least $s>0$ at the point $P$
and does~not contain components of $\bar{D}\cdot T$ that pass
through $P$, where $s$ is a natural number. Then
$$
s20\big(n/18-\mu/20\big)=\bar{D}\cdot H\cdot T\geqslant \mathrm{mult}_{P}\big(\bar{D}\big)s>\big(10n/9-\mu\big) s,%
$$
which implies that $\mu<n/9$, but $\mu>n/10$.

We see that $P\in C$. Then $\bar{M}$ contains $C$ and $L$. Put
$$
\bar{D}\big\vert_{\bar{M}}=m_{1}L+m_{2}C+\Delta,
$$
where $m_{1}$ and $m_{2}$ are a non-negative integers, and
$\Delta$ is an effective cycle, whose support does not contain $L$
and $C$. Then $m_{2}\leqslant n$ by Remark~\ref{remark:curves},
because $\alpha^{*}(-K_{X})\cdot C=1/9$.

The surface $\bar{M}$ is smooth at $P$.  So it follows from
Theorem~7.5 in \cite{Ko97} that the log pair
$$
\Big(\bar{M},\ \frac{1}{n}\bar{D}\big\vert_{\bar{M}}+\big(\mu/n-1/9\big)E\big\vert_{\bar{M}}\Big)%
$$
is not log canonical in a neighborhood of the point $P$, but
$E\vert_{\bar{M}}=L+Z$, where $Z$ is an irreducible curve that
does not pass through the point $P$. Therefore, the singularities
of the log pair
$$
\Big(\bar{M},\ \big(m_{1}/n+\mu/n-1/9\big)L+C+\frac{1}{n}\Delta\Big)%
$$
are not log canonical at the point $P$. So it follows from
Theorem~7.5 in \cite{Ko97}~that
$$
n/9-\mu-m_{1}+m_{2}=\Delta\cdot C\geqslant
\mathrm{mult}_{P}\big(\Delta\vert_{C}\big)>n-m_{1}-\mu+n/9,
$$
because $C^{2}=-1$ and $C\cdot L=1$ on $\bar{M}$. Thus, we have
$m_{2}>n$, which is a contradiction.
\end{proof}

Suppose that $\gimel=13$. Then $r=5$ by Lemma~\ref{lemma:a-not-1}.
The~base locus of the pencil $|-K_{U}|$ consists of two curves
$\bar{C}$ and $\bar{L}$ such that $\bar{C}=E\vert_{T}$, and
$\alpha(\bar{L})$~is~the~base curve of $|-K_{X}|$. Then
$$
\bar{C}^{2}=\bar{L}^{2}=-5/6,\ \bar{L}\cdot\bar{C}=1
$$
on the surface $T$. Put
$\bar{D}\vert_{T}=\bar{m}_{L}\bar{L}+\bar{m}_{C}\bar{C}+\Upsilon$,
where $\bar{m}_{L}$ and $\bar{m}_{C}$ are non-negative integers,
and $\Upsilon$ is an effective cycle, whose support does not
contain $\bar{L}$ and $\bar{C}$. Then
$$
11n/5-11\mu/6=\big(6L+5C\big)\cdot \Big(\bar{m}_{L}\bar{L}+\bar{m}_{C}\bar{C}+\Upsilon\Big)=11\bar{m}_{C}/6+\big(6L+5C\big)\cdot\Upsilon\geqslant 11\bar{m}_{C}/6,%
$$
which implies that $\bar{m}_{C}\leqslant 6n/5-\mu$. Thus, we have
 $\bar{m}_{C}<n$, because $\mu>n/5$.

Suppose that $P\not\in\bar{L}$. Then it follows from Theorem~7.5
in \cite{Ko97} that the log pair
$$
\Big(S,\ \bar{C}+\frac{\bar{m}_{L}}{n}L+\frac{1}{n}\Upsilon\Big)%
$$
is not log canonical in the neighborhood of the point $P$, because
$\bar{m}_{C}+\mu-n/5\leqslant n$, which~implies that the
inequality $\mathrm{mult}_{P}(\Upsilon\vert_{\bar{C}})>n$ holds by
Theorem~7.5 in \cite{Ko97}. Hence, we have
$$
5\mu/6+5\bar{m}_{C}/6\geqslant5\mu/6-\bar{m}_{L}+5\bar{m}_{C}/6=\Upsilon\cdot\bar{C}>n,
$$
which is impossible, because $\bar{m}_{C}\leqslant 6/5-\mu$. Thus,
we see that $P=\bar{L}\cap\bar{C}$.

Put $\bar{D}\vert_{\bar{M}}=m\bar{L}+\Omega$, where $m$ is a
non-negative integer, and $\Omega$ is an effective cycle, whose
support does not contain $\bar{L}$. Then $L^2=1/6$ on the surface
$\bar{M}$. But $m\leqslant n$ by Remark~\ref{remark:curves}.

Arguing as in the case $P\not\in\bar{L}$ we see that
$\mathrm{mult}_{P}(\Omega\vert_{\bar{L}})>n$ by Theorem~7.5 in
\cite{Ko97}, and
$$
11n/30-\mu=\bar{D}\cdot\bar{L}=m/6+\Omega\cdot\bar{L}>m/6+n,
$$
which implies that $m<0$.  So we have a contradiction, which
completes the proof of Lemma~\ref{lemma:quadratic-involutions}.

\section{Direct products.}
\label{section:conic-bundles}

Let $X$ be an arbitrary Fano variety with terminal
$\mathbb{Q}$-factorial singularities of Picard rank one, and
$\Gamma$ be a subgroup of the group $\mathrm{Bir}(X)$.

\begin{definition}
\label{definition:untwist} The subgroup
$\Gamma\subset\mathrm{Bir}(X)$ untwists all maximal singularities
if for every linear system $\mathcal{M}$ on the variety $X$ that
has no fixed components there~is $\xi\in\Gamma$ such that
the~singularities of the log pair $(X, \lambda\xi(\mathcal{M}))$
are canonical, where $\lambda\in\mathbb{Q}$ such~that
$K_{X}+\lambda\,\xi(\mathcal{M})\equiv 0$.
\end{definition}

It is well known that the group $\mathrm{Bir}(X)$ is generated by
the subgroups $\Gamma$ and $\mathrm{Aut}(X)$ in the~case when the
subgroup $\Gamma$ untwists all maximal singularities (see
\cite{Ch05umn}).

\begin{definition}
\label{definition:birational-rigidity} The variety $X$ is
birationally superrigid\footnote{There are several similar
definitions of birational superrigidity and birational
rigidity~(see~\cite{CPR},~\cite{Ch05umn}).} (rigid, respectively)
if the~trivial~subgroup (the whole group $\mathrm{Bir}(X)$,
respectively) untwists all maximal singularities.
\end{definition}

The birational rigidity of $X$ implies that there is no dominant
rational map $\rho\colon X\dasharrow Y$ such that
$\mathrm{dim}(Y)\geqslant 1$, and sufficiently general fiber of
the map $\rho$ is rationally connected (see \cite{Ch05umn}).

\begin{example}
\label{example:Pukhlikov} It follows from \cite{Pu04d} that the
variety $X$ is birationally superrigid and $\mathrm{lct}(X)=1$ in
the case when $X$ is one of the following smooth Fano varieties:
\begin{itemize}
\item a general hypersurface in $\mathbb{P}^{r}$ of degree $r\geqslant 6$;%
\item a general hypersurface in $\mathbb{P}(1^{m+1}, m)$ of degree $2m\geqslant 6$.%
\end{itemize}
\end{example}

\begin{definition}
\label{definition:universally-untwist} The subgroup $\Gamma$
universally untwists all maximal singularities if for every
variety $U$, and every linear system $\mathcal{M}$ on the variety
$X\times U$ that does not have~fixed components, there is a
birational automorphism $\xi\in\Gamma$ such that the log pair
$$
\Big(F,\ \lambda\xi\big(\mathcal{M}\big)\vert_{F}\Big)
$$
has at most canonical singularities, where $F$~is~a sufficiently
general fiber of the natural projection~$X\times U\to U$, and
$\lambda$ is a positive rational number such~that
$K_{F}+\lambda\,\xi(\mathcal{M})\vert_{F}\equiv 0$.
\end{definition}

Let $X_{1},\ldots,X_{r}$ be Fano varieties of Picard rank one with
terminal $\mathbb{Q}$-fac\-to\-ri\-al singularities.~Put
$$
U_{i}=X_{1}\times\cdots\times X_{i-1}\times \widehat{X_{i}}\times X_{i+1}\times\cdots\times X_{r}%
$$
and $V=X_{1}\times\cdots\times X_{r}$. Let $\pi_{i}\colon V\to
U_{i}$ be a natural projection.

For every $i\in\{1,\ldots,r\}$,~suppose that
$\mathrm{lct}(X_{i})\geqslant 1$,~and~there is a~subgroup
$\Gamma_{i}\subset\mathrm{Bir}(X_{i})$ that universally~untwists
all maximal singularities. Then the following result
holds\footnote{The assertion of Theorem~\ref{theorem:Cheltsov} is
proved in \cite{Pu04d} for smooth birationally superrigid Fano
varieties.}.

\begin{theorem}
\label{theorem:Cheltsov} The variety $X_{1}\times\cdots\times
X_{r}$ is non-rational, and
$$
\mathrm{Bir}\Big(X_{1}\times\cdots\times X_{r}\Big)=\Big<\prod_{i=1}^{r}\Gamma_{i},\ \mathrm{Aut}\Big(X_{1}\times\cdots\times X_{r}\Big)\Big>,%
$$%
for any dominant rational map $\rho\colon X_{1}\times\cdots\times
X_{r}\dasharrow Y$ whose general fiber is rationally connected,
there is a commutative diagram
$$
\xymatrix{
X_{1}\times\cdots\times X_{r}\ar@{->}[d]_{\pi}\ar@{-->}[rr]^{\sigma}&&X_{1}\times\cdots\times X_{r}\ar@{-->}[d]^{\rho}\\
X_{i_{1}}\times\cdots\times X_{i_{k}}\ar@{-->}[rr]_{\xi}&&Y,}%
$$
where $\xi$ and $\sigma$ are
birational maps, and $\pi$ is a  projection for some $\{i_{1},\ldots,i_{k}\}\subseteq\{1,\ldots,r\}$.%
\end{theorem}

It is well known that Theorem~\ref{theorem:Cheltsov} is implied by
the following technical result (see \cite{Pu04d}).

\begin{proposition}
\label{proposition:Cheltsov} For every linear system $\mathcal{M}$
on the variety $V$ such that
\begin{itemize}
\item the linear system $\mathcal{M}$ does not have fixed components,%
\item the linear system $\mathcal{M}$ does not lie in the fibers of the projections $\pi_{1},\ldots,\pi_{r}$,%
\end{itemize}
there are $k\in\{1,\ldots,r\}$, birational map
$\xi\in\prod_{i=1}^{r}\Gamma_{i}$ and a positive rational number
$\mu$ such~that
\begin{itemize}
\item the inequality $\kappa(V, \mu\xi(\mathcal{M}))\geqslant 0$ holds\,\footnote{The number $\kappa(V,\ \mu\xi(\mathcal{M}))$ is a Kodaira dimension of the movable log pair $(V,\ \mu\xi(\mathcal{M}))$ (see \cite{Ch05umn}).},%
\item the equivalence $K_{V}+\mu\,\xi(\mathcal{M})\equiv \pi_{k}^{*}\big(D\big)$ holds for some nef $\mathbb{Q}$-divisor $D$ on $U_{k}$.%
\end{itemize}
\end{proposition}

\begin{proof}
Let $F_{i}$ be a sufficiently general fiber of $\pi_{i}$. The
subgroups $\Gamma_{1},\ldots,\Gamma_{r}$ universally~untwist all
maximal singularities for every $i=1,\ldots,r$. So there~is
$\xi\in\prod_{i=1}^{r}\Gamma_{i}$ such that the log pairs
$$
\Big(F_{1},\ \mu_{1}\,\xi\big(\mathcal{M}\big)\big\vert_{F_{1}}\Big),\ \ldots,\ \Big(F_{r},\ \mu_{r}\,\xi\big(\mathcal{M}\big)\big\vert_{F_{r}}\Big)%
$$
are canonical, where $\mu_{i}$ is a rational number such that
$$
K_{V}+\mu_{i}\,\xi\big(\mathcal{M}\big)\equiv \pi_{i}^{*}\big(D_{i}\big),%
$$
where $D_{i}$ is a $\mathbb{Q}$-divisor on $U_{i}$. Then there is
$m\in\{1,\ldots,r\}$ such that $D_{m}$ is nef.

Now arguing as in the proof of Theorem~1 in \cite{Pu04d}, we see
that $\kappa(V, \mu_{k}\xi(\mathcal{M}))\geqslant 0$.
\end{proof}

Let $X$ be a general quasismooth hypersurface in
$\mathbb{P}(1,a_{1},a_{2},a_{3},a_{4})$ of degree
$\sum_{i=1}^{4}a_{i}$ with terminal singularities, where
$a_{1}\leqslant a_{2}\leqslant a_{3}\leqslant a_{4}$. Then $X$ is
a Fano threefold, whose divisor class group is generated by
$-K_{X}$. The possible values of $(a_{1},a_{2},a_{3},a_{4})$ are
given in Table~5 in \cite{IF00}.

There are finitely many non-biregular birational involutions
$\tau_{1},\ldots,\tau_{k}\in\mathrm{Bir}(X)$
explicitly~constructed in \cite{CPR}~such~that the following
result holds (see \cite{CPR}).

\begin{theorem}
\label{theorem:CPR} The subgroup
$\langle\tau_{1},\ldots,\tau_{k}\rangle$ universally untwists all
maximal singularities.
\end{theorem}

Hence, the following two examples follow from \cite{MaMon64},
\cite{ChPa05}, \cite{Pu04d} and Theorems~\ref{theorem:main},
\ref{theorem:Cheltsov}, \ref{theorem:CPR}.

\begin{example}
\label{example:41} Let $X$ be a general hypersurface in
$\mathbb{P}(1,1,4,5,10)$ of degree~$20$, and $U$~be~a~general
hypersurface in $\mathbb{P}(1^{n+1}, n)$ of degree $2n\geqslant
6$. Then $\mathrm{Bir}(X\times
U)\cong(\mathbb{Z}_{2}\ast\mathbb{Z}_{2})\oplus\mathbb{Z}_{2}$.
\end{example}

\begin{example}
\label{example:9} Let $X$ be a general hypersurface in
$\mathbb{P}(1,1,2,3,3)$ of degree~$9$, and $U$~be~a~general
hypersurface in $\mathbb{P}^{r}$ of degree $r\geqslant 6$. Then
$$
\mathrm{Bir}\Big(X\times U\Big)\cong\Big<a,b,c\ \big\vert\ a^2=b^2=c^2=\big(abc\big)^{2}=1\Big>.%
$$
\end{example}

It follows from \cite{CPR} that
$\mathrm{Aut}(X)\ne\mathrm{Bir}(X)$ for exactly $45$ values of
$(a_{1},a_{2},a_{3},a_{4})$.

\end{document}